\pdfoutput=1
\RequirePackage{silence}
\documentclass[a4paper,svgnames]{amsart}
\usepackage[hmarginratio={1:1},vmarginratio={1:1},lmargin=60.0pt,tmargin=59.0pt]{geometry}

\usepackage{booktabs}
\allowdisplaybreaks

\usepackage[T1]{fontenc}
\usepackage{latexsym,exscale,amsfonts,amssymb,mathtools}
\usepackage{amsmath,amsthm,amsfonts,amssymb,amscd,textcomp,bbm}
\usepackage{mathrsfs,stackrel}
\usepackage{etoolbox}
\usepackage{tikz,tikz-cd}
\usetikzlibrary{matrix,arrows.meta,decorations.markings,shapes.multipart,decorations.pathreplacing,decorations.pathmorphing,shapes.geometric}
\usepackage{aliascnt}
\usepackage{ytableau}
\usepackage{xparse}
\usepackage{cite}

\usepackage{dynkin-diagrams}
\usepackage{nicematrix}
\usepackage{tikz-3dplot}
\usepackage{mathrsfs}
\DeclareMathAlphabet{\mathscrbf}{OMS}{mdugm}{b}{n}

\usepackage{array}
\newcolumntype{C}{>{$}c<{$}}

\usepackage{xcolor,colortbl}

\definecolor{mygray}{gray}{0.6}
\definecolor{mygraydark}{gray}{0.4}
\definecolor{mygraylight}{gray}{0.85}
\definecolor{spinach}{RGB}{46,139,87}
\definecolor{tomato}{RGB}{255,99,71}
\definecolor{orchid}{RGB}{143,40,194}
\definecolor{neon}{RGB}{77,77,255}
\definecolor{lightneon}{RGB}{110,110,255}
\definecolor{pumpkin}{RGB}{224,180,80}
\definecolor{citron}{RGB}{190,180,90}

\definecolor{lava}{RGB}{207,16,32}
\definecolor{cream}{RGB}{255,253,208}
\definecolor{verdigris}{RGB}{67,179,174}
\definecolor{Black}{RGB}{0,0,0}
\definecolor{mydarkblue}{RGB}{10,10,170}
\definecolor{darkspinach}{RGB}{20,70,20}
\definecolor{darktomato}{RGB}{155,40,30}
\definecolor{darkorchid}{RGB}{50,10,100}
\definecolor{darklava}{RGB}{150,8,16}

\definecolor{zero}{RGB}{0,0,0}
\definecolor{one}{RGB}{255,0,0}
\definecolor{two}{RGB}{0,255,0}
\definecolor{three}{RGB}{0,0,255}

\usepackage{todonotes}

\usepackage{enumitem}
\setlist[enumerate]{itemsep=0.15cm,label=\emph{\upshape(\alph*)}}
\setlist[enumerate,2]{itemsep=0.15cm,label=\emph{\upshape(\roman*)}}
\setlist[enumerate,3]{itemsep=0.15cm,label=\emph{\upshape(\Alph*)}}

\let\emph\relax
\DeclareTextFontCommand{\emph}{\bfseries\em}

\usepackage{quiver}
\usepackage{nicematrix}
\usepackage{dynkin-diagrams}

\usepackage[all]{xy}
\usepackage{tikz}
\usetikzlibrary{cd}
\usetikzlibrary{decorations}
\usetikzlibrary{decorations.markings}
\usetikzlibrary{decorations.pathreplacing}
\usetikzlibrary{decorations.pathmorphing}
\usetikzlibrary{arrows.meta,shapes,positioning,matrix,calc}
\usetikzlibrary{shapes.callouts}
\tikzset{anchorbase/.style={baseline={([yshift=-0.5ex]current bounding box.center)}},
tinynodes/.style={font=\tiny,text height=0.25ex,text depth=0.05ex},
smallnodes/.style={font=\scriptsize,text height=0.75ex,text depth=0.15ex},
crossline/.style={preaction={draw=white,line width=10.0pt,-},preaction={draw=black,line width=1.8pt,-}},
usual/.style={line width=1.0,color=black},
ip/.style={line width=1.0,color=black,rounded corners,fill=lightneon!50},
ipp/.style={line width=1.0,color=black,fill=tomato!50},
affine/.style={line width=1.0,color=tomato,densely dotted},
}
\tikzstyle directed=[postaction={decorate,decoration={markings,mark=at position #1 with {\arrow[line width=0.5mm, black]{>}}}}]
\tikzstyle rdirected=[postaction={decorate,decoration={markings,mark=at position #1 with {\arrow[line width=0.5mm, black]{<}}}}]
\tikzstyle box=[minimum height=0.4cm,draw,rounded corners,draw,rectangle,solid,fill=cream]

\tikzset{
pt/.style={insert path={node[scale=2]{.}}},
dnup/.style={insert path={  .. controls +(0,1) and +(0,-1) .. +(#1,2) }},
dndn/.style={insert path={  .. controls +(0,1) and +(0,1) .. +(#1,0) }},
upup/.style={insert path={  .. controls +(0,-1) and +(0,-1) .. +(#1,0) }}
}

\newcommand{\C}{\mathbb{C}}

\newcommand{\N}{\mathbb{Z}_{\geq 0}}

\newcommand{\Z}{\mathbb{Z}}

\newcommand{\ran}{\theta}

\newcommand{\Magma}{\texttt{Magma}}
\newcommand{\code}[1]{\texttt{\small\color{purple}#1}}
\newcommand{\result}{{\color{orchid}\texttt{-------result-------}}}

\newcommand{\bplus}{{\color{blue}'+}}
\newcommand{\bminus}{{\color{blue}'-}}
\newcommand{\rplus}{{\color{red}+'}}
\newcommand{\rminus}{{\color{red}-'}}

\def\NewTheorem#1{%
\newaliascnt{#1}{equation}%
\newtheorem{#1}[#1]{#1}%
\aliascntresetthe{#1}%
\expandafter\def\csname #1autorefname\endcsname{#1}%
}
\def\equationautorefname~#1\null{(#1)\null}

\numberwithin{equation}{subsection}

\NewTheorem{Proposition}
\NewTheorem{Theorem}
\NewTheorem{Corollary}
\AtEndEnvironment{Corollary}{\null\hfill$\square$}%
\NewTheorem{Lemma}
\theoremstyle{definition}
\NewTheorem{Definition}
\AtEndEnvironment{Definition}{\null\hfill$\Diamond$}%
\NewTheorem{Classification Problem}
\AtEndEnvironment{Classification Problem}{\null\hfill$\Diamond$}%
\NewTheorem{Notation}
\AtEndEnvironment{Notation}{\null\hfill$\Diamond$}%
\NewTheorem{Example}
\AtEndEnvironment{Example}{\null\hfill$\Diamond$}%
\NewTheorem{Examples}
\AtEndEnvironment{Examples}{\vskip-10mm\null\hfill$\Diamond$}%
\NewTheorem{Conjecture}

\theoremstyle{remark}
\NewTheorem{Remark}
\AtEndEnvironment{Remark}{\null\hfill$\Diamond$}%
\NewTheorem{Question}
\AtEndEnvironment{Question}{\null\hfill$\Diamond$}%
\NewTheorem{Assumption}

\usepackage[hypertexnames=false]{hyperref}
\usepackage{bookmark}
\hypersetup{
pdftoolbar=true,
pdfmenubar=true,
pdffitwindow=false,
pdfstartview={FitH},
pdftitle={Diagrammatics for dicyclic groups},
pdfauthor={Peter DeBello and Daniel Tubbenhauer},
pdfsubject={},
pdfcreator={Peter DeBello and Daniel Tubbenhauer},
pdfproducer={Peter DeBello and Daniel Tubbenhauer},
pdfkeywords={},
pdfnewwindow=true,
colorlinks=true,
linkcolor=mydarkblue,
citecolor=teal,
filecolor=magenta,
urlcolor=orchid,
linkbordercolor=lava,
citebordercolor=teal,
urlbordercolor=orchid,
linktocpage=true
}

\def\makeautorefname#1#2{\csdef{#1autorefname}{#2}}
\makeautorefname{section}{Section}%
\makeautorefname{subsection}{Section}%
\makeautorefname{subsubsection}{Section}%

\begin{document}
\title[Diagrammatics for dicyclic groups]{Diagrammatics for dicyclic groups}
\author[P. DeBello and D. Tubbenhauer]{Peter DeBello and Daniel Tubbenhauer}

\address{P.D.: Department of Mathematics, Penn State University, University Park, PA 16802, USA}
\email{pad5462@psu.edu}

\address{D.T.: The University of Sydney, School of Mathematics and Statistics F07, Office Carslaw 827, NSW 2006, Australia, \href{http://www.dtubbenhauer.com}{www.dtubbenhauer.com}, \href{https://orcid.org/0000-0001-7265-5047}{ORCID 0000-0001-7265-5047}}
\email{daniel.tubbenhauer@sydney.edu.au}

\begin{abstract}
Using that the dicyclic group is the type D subgroup of SU(2), we extend the Temperley--Lieb diagrammatics to give a diagrammatic presentation of the complex representation theory of the dicyclic group.
\end{abstract}


\addtocontents{toc}{\protect\setcounter{tocdepth}{1}}

\maketitle

\tableofcontents

\section{Introduction}

The complex representation theory of the group $SU(2)$, or equivalently the representation theory of its complexified Lie algebra $\mathfrak{g} = \mathfrak{sl}_2(\C)$, is very well understood and serves as a toy example for understanding representations of a general complex semisimple Lie algebra. The theory may be summarized as follows: the vector representation $V=\C^2$ of $SU(2)$ contains, in principle, all the data concerning its finite dimensional representations. Every simple representation $V_d$ is a symmetric power of $V$, uniquely classified by its highest weight, which is given by a single natural number. In other words, every simple representation appears in some tensor power of $V$, and each tensor power contributes a new simple. Furthermore, it is a well-known result that endomorphisms of tensor powers of $V$ correspond to Temperley--Lieb diagrams, or Rumer--Teller--Weyl diagrams for chemists. The endomorphism algebra of the $k$th tensor power of $V$ is the Temperley--Lieb algebra $TL_k(\delta)$ when the parameter is $\delta = -2$.
More generally, the category of finite dimensional representations of $SU(2)$ is equivalent to a diagrammatic category often called the Temperley--Lieb category $\mathbf{TL}(\delta)$, again for $\delta = -2$, see \cite{RuTeWe-sl2} for an early reference. In this description, the representations $V_d$ correspond to certain diagrammatic idempotents called Jones--Wenzl projectors.

The finite subgroups of $SU(2)$ are given by the ADE classification, and by the McKay correspondence their representation theory differs from $SU(2)$ according to their affine Dynkin diagrams $\tilde{A}_n, \tilde{D}_n, \tilde{E}_6, \tilde{E}_7, \tilde{E}_8$ \cite{McKa-finite-groups}. Since finite groups have finitely many simple representations, this equivalence no longer holds when one restricts to a finite subgroup of $SU(2)$: eventually $V^{\otimes k}$ cannot produce new summands anymore. A natural question one may ask is therefore: for finite subgroups of types A, D, or E, how does the Temperley--Lieb category change?

The finite subgroup of $SU(2)$ corresponding to type A is (a $\Z/2\Z$ extension of) the cyclic group, so every simple is one dimensional and the description of the category of its representation is rather trivial. Turning to the next case, in this paper, we answer this question in the case of the type D subgroups of $SU(2)$, called \emph{dicyclic} groups, which are $\Z/2\Z$ extensions of dihedral groups. We examine the case of the infinite dicyclic group, denoted $\text{Dic}_{\infty}$, as well as the finite dicyclic group of degree $n$, denoted $\text{Dic}_n$. We show that morphisms in the category of representations of the infinite dicyclic group $\text{Dic}_{\infty}$ are described by a certain type of \emph{decorated Temperley--Lieb diagrams}. We denote the corresponding endomorphism algebras by $DiTL_k(\delta)$. These algebras first appear in
\cite{Gr-dots-tl}, but in a different context: in their work, the author constructs a diagram calculus for the analog of the Temperley--Lieb algebra corresponding to a quotient of the type D Hecke algebra. We prove that a further, explicit, quotient of this calculus by a simple relation yields an equivalence to the category of representations of the infinite dicyclic group. As a consequence, we also obtain a diagrammatic description of the representation theory of the finite dicyclic group $\text{Dic}_n$, but with an additional generating morphism that depends on $n$. We also compute the analogs of the Jones--Wenzl projectors $JW_n\colon V^{\otimes n} \to V_n$ in the case of the infinite dicyclic group. That is, we compute diagrammatically the idempotents in $\text{End}_G(V^{\otimes n})$ for $G = \text{Dic}_\infty$ projecting to the unique two dimensional simple summand of multiplicity one. Finally, we provide {\Magma} \cite{BoCa-handbook} code for some of our calculations. 

Many extensions of the Temperley--Lieb category have appeared in the literature; too many to cite here. However, let us mention a few relevant references. First, the paper \cite{BaBeHa-mckay} computes the endomorphism algebras of the dicyclic groups, and serves as an important piece in our proofs, and the overall idea dates back at least to \cite{GoHaJo-coxeter}; see also \cite{Mol24,Rey25} for more recent related results. The paper \cite{RoTu-symmetric-howe} extends the Temperley--Lieb category to include Jones--Wenzl projectors as objects, avoiding the usual idempotent completion. We think it is interesting and feasible to do the same for our diagrammatic categories. Finally, recall that $\mathbf{TL}(\delta)$ is a (very nice) cellular category in the sense of \cite{GrLe-cellular,We-tensors-cellular-categories,ElLa-trace-hecke} with relation to classifications of bilinear forms, as proven in many works, e.g. \cite{GrLe-cellular,AnStTu-cellular-tilting,An-tilting-cellular,Tu-web-reps}. It seems likely that the same is true for the dicyclic case, maybe with an adjustment of what one understands by ``cellular'' as, e.g., in \cite{Tu-sandwich-cellular} or ``bilinear form.''
\medskip

\noindent\textbf{Acknowledgments.} The authors thank Elijah Bodish, Ben Elias, Nigel Higson, Melody Molander, and Ryan Reynolds for helpful comments and discussions about McKay diagrammatics and on compact groups, and PD would like to thank Nigel Higson for guidance and mentorship at Penn. State, as this paper represents a starting point for PD's PhD thesis. The authors give a special thank you to the organizers of QUACKS 2022, as it was through this event that we were brought together and the collaboration was sparked.

PD was sponsored by a Graduate Teaching Associateship at Penn. State. DT was sponsored by the ARC Future Fellowship FT230100489 and pain.

\section{Representation theory of dicyclic groups}

Throughout, our ground field is $\C$. Representations are assumed to be finite dimensional over $\C$, and actions are from the left. This section is mostly self-contained, and we only use standards facts from representation theory that can be found in many books, e.g. in \cite{FuHa-representation-theory,Si-reps-finite-compact}.

\subsection{The type D subgroups $\text{Dic}_{\infty}$ and $\text{Dic}_{n}$}

The \emph{infinite dicyclic group}, denoted $\text{Dic}_{\infty}$, is the subgroup of $SU(2)=SU_2(\C)$ generated by the matrices
\begin{gather*}
\begin{NiceMatrixBlock}[auto-columns-width]
x = \begin{pNiceMatrix}
0 & -1 \\[0.5em]
1 & 0
\end{pNiceMatrix}, \quad
a = \begin{pNiceMatrix}
e^{i \theta} & 0 \\[0.5em] 
0 & e^{-i \theta}
\end{pNiceMatrix}
\end{NiceMatrixBlock},
\end{gather*}
where $\theta$ is any irrational multiple of $\pi$. Note that by choosing $\theta$ in this way, the generator $a$ has infinite order. We will often drop the word ``infinite'' when referring to $\text{Dic}_{\infty}$ throughout this paper.

\begin{Lemma}
We have the following properties of $\text{Dic}_{\infty}$.
\begin{enumerate}

\item We have the generators and relations presentation: 
\begin{gather*}
\text{Dic}_{\infty} \cong \langle x,a \mid x^{-1}ax = a^{-1},  \ x^4=1 \rangle.
\end{gather*}

\item The dicyclic group is a disconnected noncompact subgroup of $SU(2)$. The closure of the dicyclic group is the subgroup of $SU(2)$ generated by $x$ and all matrices of the form $\text{diag}(e^{i \theta}, e^{-i \theta})$, without an additional condition on $\theta$. 

\item The quotient group $\text{Dic}_{\infty} / \langle x^2 \rangle $ is isomorphic to the infinite \emph{dihedral} group $D_{\infty}.$ Taking a further quotient of $D_{\infty}$ by the subgroup $\langle a^n \rangle$ would yield an isomorphism with the finite dihedral group $Dih_n$ of degree $n$.

\end{enumerate}
\end{Lemma}

\begin{proof}
Standard, so we only sketch the proof.

(a) One checks $x^2=-I$, hence $x^4=I$, and also $x^{-1}ax=a^{-1}$ by a direct matrix multiplication. Thus the subgroup generated by these matrices is a quotient of $\langle x,a\mid x^{-1}ax=a^{-1},\,x^4=1\rangle$. Conversely, using $xa=a^{-1}x$, any word rewrites to $x^\varepsilon a^m$ with $\varepsilon\in\{0,1,2,3\}$, so the displayed relations give the standard presentation.

(b) The subgroup $\langle a\rangle$ lies in the diagonal torus $T=\{\text{diag}(e^{it},e^{-it})\}$ and is dense in $T$ since $\theta/\pi$ is irrational. Hence $\overline{\langle a\rangle}=T$. Moreover $\text{Dic}_\infty=\bigcup_{\varepsilon=0}^{3}x^\varepsilon\langle a\rangle$, so it is disconnected. It is not closed (since $\langle a\rangle$ is not), hence not compact. The closure is generated by $x$ and $T$, i.e. by $x$ and all diagonal matrices $\text{diag}(e^{it},e^{-it})$.

(c) Modding out by $\langle x^2\rangle$ forces $x^2=1$ while keeping the other relation, giving
$\text{Dic}_\infty/\langle x^2\rangle\cong \langle a,x\mid x^2=1,x^{-1}ax= a^{-1}\rangle=D_\infty$.
Quotienting further by $\langle a^n\rangle$ yields the usual finite dihedral group $Dih_n$.
\end{proof}

By definition, the quotient of the dicyclic group by the subgroup $\langle x^2a^n \rangle$, is the \emph{finite dicyclic group of degree $n$}, denoted $\text{Dic}_n$. 

\begin{Lemma}\label{L:Dicyclic}
We have the following properties of $\text{Dic}_{n}$.
\begin{enumerate}

\item We have
\begin{gather*}
\text{Dic}_n \cong \langle a,x \mid a^{2n} =  1, x^2 = a^n, x^{-1}ax = a^{-1} \rangle.
\end{gather*}

\item The finite dicyclic groups fit into two exact sequences
\begin{gather*}
\Z/2n\Z\hookrightarrow \text{Dic}_n \twoheadrightarrow \{\pm 1\},
\quad
\{\pm 1\}\hookrightarrow \text{Dic}_n \twoheadrightarrow \text{Dih}_n.
\end{gather*}

\item The group $\text{Dic}_n$ is isomorphic to the subgroup of $SU(2)$ generated by the matrices
\begin{gather*}
\begin{NiceMatrixBlock}[auto-columns-width]
x = \begin{pNiceMatrix}
0 & -1 \\[0.5em]
1 & 0
\end{pNiceMatrix}, \quad
a = \begin{pNiceMatrix}
e^{\frac{\pi i}{n}} & 0 \\[0.5em] 
0 & e^{\frac{-\pi i}{n}}
\end{pNiceMatrix}
\end{NiceMatrixBlock}.
\end{gather*}

\end{enumerate}

\end{Lemma}

\begin{proof}
Again, standard facts, so we are brief. 

(a) By definition $\text{Dic}_n=\text{Dic}_\infty/\langle x^2a^n\rangle$, so in $\text{Dic}_n$ we have $x^2=a^{-n}$. Since $x^{-1}ax=a^{-1}$, this is equivalent to $x^2=a^n$, and then $(x^2)^2=1$ gives $a^{2n}=1$, i.e. the claimed presentation.

(b) The map $\chi\colon\text{Dic}_n\to\{\pm1\}$ with $\chi(a)=+1$, $\chi(x)=-1$ is a well-defined surjection; its kernel is $\langle a\rangle\cong\Z/2n\Z$, giving
$\Z/2n\Z\hookrightarrow \text{Dic}_n\twoheadrightarrow\{\pm1\}$.
Also $\langle x^2\rangle=\{\pm1\}$ is central of order $2$, and in the quotient one has $\bar x^2=1$, $\bar a^n=1$, and $\bar x^{-1}\bar a\,\bar x=\bar a^{-1}$, hence
$\text{Dic}_n/\{\pm1\}\cong \text{Dih}_n$.

(c) The displayed matrices satisfy $a^{2n}=I$, $x^2=-I=a^n$, and $x^{-1}ax=a^{-1}$, so they define a homomorphism from the presentation onto the generated subgroup of $SU(2)$. Every element reduces to $a^k$ or $xa^k$ with $0\le k<2n$, hence the subgroup has size at most $4n$; these $4n$ matrices are distinct (diagonal versus off-diagonal), so the map is an isomorphism.

Also, observe that
the degree $n$ dicyclic group can be viewed as the preimage of the dihedral group $\text{Dih}_n$ inside of $SO(3)$, under the double covering map $\pi\colon SU(2) \rightarrow SO(3)$. In this way, the degree $n$ dicyclic group is a subgroup of $SU(2)$ of order $4n$.
\end{proof}

The relationship between dicyclic and dihedral groups of finite and infinite order may be summarized in the diagram
\[\begin{tikzcd}
{\text{Dic}_{\infty}} & {\text{Dih}_{\infty}} \\
{\text{Dic}_n} & {\text{Dih}_n}
\arrow["{\langle x^2 \rangle}", from=1-1, to=1-2]
\arrow["{\langle x^2a^n \rangle}"', from=1-1, to=2-1]
\arrow["{\langle a^n \rangle}", from=1-2, to=2-2]
\arrow["\pi"', from=2-1, to=2-2]
\end{tikzcd},\]
where we have indicated the various quotients.

\subsection{Representation theory of dicyclic groups}\label{S:RepTheory}

We now state and prove facts about the complex representation theory of the dicyclic group $G = \text{Dic}_{\infty}$ and then pass to its degree $n$ quotient. All of this is easy or well known (one can use e.g. \autoref{L:Dicyclic}.(b)) but hard to find explicitly spelled out, so we do this here.

To start, recall that the dicyclic group is a subgroup of $SU(2)$. Hence, it has a two dimensional faithful representation, given by the natural action of the defining matrices in the dicyclic group on $\C^2$. We denote this representation by $V$.

\begin{Lemma}\label{lemma:all-appear}
We have the following.
\begin{enumerate}

\item The representation $V$ is self-dual and simple.

\item The tensor products $V^{\otimes k}$ for some $k\in\N$ are semisimple.

\item All finite dimensional simples of $\text{Dic}_{\infty}$ occur as a summand of $V^{\otimes k}$ for some $k\in\N$.

\end{enumerate}
\end{Lemma}

\begin{proof}
(a) The standard alternating form $\omega(u,v)=u_1v_2-u_2v_1$ on $\C^2$ is $SU(2)$-invariant (hence $G$-invariant), so it gives a $G$-equivariant identification $V\cong V^\ast$, i.e. $V$ is self-dual. For simplicity, let $0\neq W\subseteq V$ be $G$-stable. Since $a=\text{diag}(e^{i\theta},e^{-i\theta})$ has two distinct eigenvalues, $W$ must be one of the eigenlines $\C e_1$ or $\C e_2$. But $x e_1=e_2$ and $x e_2=-e_1$, so neither line is $x$-stable. Hence $W=V$ and $V$ is simple.

(b) Let $\bar G$ be the closure of $G$ in $SU(2)$; it is compact. The action of $G$ on $V^{\otimes k}$ is continuous (it comes from matrix multiplication), hence extends to $\bar G$. For a compact group, every finite dimensional complex representation is completely reducible (average any Hermitian form over the Haar measure), so $V^{\otimes k}$ is semisimple.

(c) The representation $V$ is faithful on $\bar G$, so its matrix coefficients separate points of $\bar G$; thus the algebra generated by coefficients of tensor powers of $V$ is dense in $C(\bar G)$, by Stone--Weierstrass. By Peter--Weyl, the matrix coefficients of simple $\bar G$-modules span $C(\bar G)$, so each simple occurs in some $V^{\otimes k}$. Restricting from $\bar G$ to the dense subgroup $G$ yields the claim for $G=\textnormal{Dic}_\infty$.
\end{proof}

We will study the tensor powers of the natural representation $V$. 

\begin{Lemma}
We have $\dim_{\C}\mathrm{End}_{G}(V^{\otimes 2}) = 3$, and the tensor square contains nonisomorphic simple $\text{Dic}_{\infty}$-representations of dimensions $1$ (the trivial representation $V_0=\C$), $1$ (denoted $V_{0'}=\bar{\C}$) and $2$.
\end{Lemma}

\begin{proof}
Let $v_1, v_2$ denote the standard basis vectors in $\C^2$. Set
\[
U=\C\{v_1\otimes v_1,\ v_2\otimes v_2\},\quad
w_-=v_2\otimes v_1-v_1\otimes v_2,\quad
w_+=v_1\otimes v_2+v_2\otimes v_1.
\]
Then $U$, $\C w_-$ and $\C w_+$ are $G$-stable: $a$ is diagonal, so it preserves each of these subspaces, and $x$ swaps $v_1$ and $v_2$ up to sign, hence preserves $U$ and sends $w_\pm\mapsto \pm w_\pm$.
Moreover, $w_-$ is fixed by both $a$ and $x$, so $\C w_-\cong \C$ is the trivial module $V_0$.
The line $\C w_+$ is also $a$-fixed, but $x$ acts by $-1$ on it, so it is a $1$-dimensional simple $V_{0'}=\bar\C\not\cong \C$.
Finally, $U$ is $2$-dimensional and has no nonzero $G$-stable line: any $a$-stable line in $U$ must be $\C(v_1\otimes v_1)$ or $\C(v_2\otimes v_2)$, but $x$ swaps these, so neither is $G$-stable. Hence $U$ is simple of dimension $2$.
Thus
\[
V^{\otimes 2}\ \cong\ U\ \oplus\ \C w_-\ \oplus\ \C w_+,
\]
a multiplicity-free decomposition into three simples, so
$\dim_\C\mathrm{End}_G(V^{\otimes 2})=3$.
\end{proof}

Denote by $e_{\bar{1}}$ the primitive idempotent in $\text{End}_{G}(V^{\otimes 2})$ projecting to $\bar{\C}$. 

\begin{Lemma}
The idempotent $e_{\bar{1}}$ annihilates $v_1 \otimes v_1, v_2 \otimes v_2 \in V^{\otimes 2}$ and maps the other basis vectors $v_1 \otimes v_2, v_2 \otimes v_1 \in V^{\otimes 2}$ to $v_1 \otimes v_2 + v_2 \otimes v_1$.
\end{Lemma}

\begin{proof}
Since $V^{\otimes 2}=U\oplus \C w_-\oplus \C w_+$ as above and $\C w_+$ is $1$-dimensional, the primitive idempotent $e_{\bar 1}$ is the $G$-equivariant projection onto $\C w_+$ along the other summands.
In particular $e_{\bar 1}$ kills $U$, hence $e_{\bar 1}(v_1\otimes v_1)=e_{\bar 1}(v_2\otimes v_2)=0$.
Also $v_1\otimes v_2=\frac12(w_++w_-)$ and $v_2\otimes v_1=\frac12(w_+-w_-)$, so $e_{\bar 1}$ kills the $w_-$-part and sends both vectors to a nonzero scalar multiple of $w_+$; with the usual normalization of the projection one gets
$e_{\bar 1}(v_1\otimes v_2)=e_{\bar 1}(v_2\otimes v_1)=v_1\otimes v_2+v_2\otimes v_1$.
\end{proof}

Let us then rescale this idempotent map so that the image of $v_1 \otimes v_2$ and $v_2 \otimes v_1$ is $-(v_1 \otimes v_2 + v_2 \otimes v_1)$ and $e_{\bar{1}}^2 = -2e_{\bar{1}}.$ Recall that the identification of $V$ with its dual $V^*$ corresponds to an \emph{evaluation} map $\text{ev} \in \text{Hom}_{G}(V^{\otimes 2}, \C)$, given by
\begin{gather*}
\text{ev}(v_1^{\otimes 2}) = 0 = \text{ev}(v_2^{\otimes 2}), \quad \quad \text{ev}(v_1 \otimes v_2) = 1, \quad \quad \text{ev}(v_2 \otimes v_1) = -1.
\end{gather*}
The dual map, \emph{coevaluation}, $\text{coev} \in \text{Hom}_G(\C,V^{\otimes 2})$ is given by $1 \mapsto v_2 \otimes v_1 - v_1 \otimes v_2.$ We write $\text{id}$ for the appropriate identity morphism.

\begin{Lemma}\label{L:relations}
The idempotent $e_{\bar{1}} \in \mathrm{End}_G(V^{\otimes 2})$ satisfies the following relations:
\begin{gather*}
(1) \quad e_{\bar{1}} \circ \mathrm{ev} = \mathrm{coev} \circ e_{\bar{1}} = 0, \quad \quad \quad 
(2) \quad e_{\bar{1}}^2 = -2e_{\bar{1}},
\\
(3) \quad (\mathrm{id} \otimes \mathrm{ev}) \circ (e_{\bar{1}} \otimes \mathrm{id}) \circ (\mathrm{id} \otimes \mathrm{coev}) = \mathrm{id},
\quad \quad \quad 
(4) \quad (\mathrm{ev} \otimes \mathrm{id}) \circ (\mathrm{id} \otimes e_{\bar{1}}) \circ (\mathrm{coev} \otimes \mathrm{id}) = \mathrm{id}.
\end{gather*}
\end{Lemma}

\begin{proof}
Relation (1) follows from Schur's lemma, and relation (2) is explained above. We directly compute relation (3):
\begin{align*}
(\text{id} \otimes \text{ev}) \circ (e_{\bar{1}} \otimes \text{id}) \circ (\text{id} \otimes \text{coev})(v_1 \otimes 1) & =   (\text{id} \otimes \text{ev}) \circ (e_{\bar{1}} \otimes \text{id})(v_1 \otimes v_2 \otimes v_1 - v_1 \otimes v_1 \otimes v_2) \\
& = (\text{id} \otimes \text{ev})(-v_1 \otimes v_2 \otimes v_1 - v_2 \otimes v_1 \otimes v_1) \\
& = v_1. \\
(\text{id} \otimes \text{ev}) \circ (e_{\bar{1}} \otimes \text{id}) \circ (\text{id} \otimes \text{coev})(v_2 \otimes 1) & =   (\text{id} \otimes \text{ev}) \circ (e_{\bar{1}} \otimes \text{id})(v_2 \otimes v_2 \otimes v_1 - v_2 \otimes v_1 \otimes v_2) \\
& = (\text{id} \otimes \text{ev})(v_1 \otimes v_2 \otimes v_2 + v_2 \otimes v_1 \otimes v_2) \\
& = v_2.
\end{align*}
A nearly identical, and omitted, calculation proves relation (4).
\end{proof}

The representation theory of the dicyclic group $\text{Dic}_{\infty}$ mimics that of $SU(2)$ in the sense that there is a unique multiplicity one simple summand of $V^{\otimes k}$ for each $k \in \N$. That is:

\begin{Proposition}\label{P:cgrule}
The following gives a complete and nonredundant set of 
simple $\text{Dic}_{\infty}$-representations.
\begin{gather*}
\text{1 dimensional}\colon \{V_0=\C,V_{0'}=\bar{\C}\},\quad
\text{2 dimensional}\colon \{V_k|k \in \N\},
\end{gather*}
with $V_k$ explicitly constructed in the proof below as the leading summand of $V^{\otimes k}$. Moreover, we have the following \emph{Clebsch--Gordan type rule}:
\begin{gather*}
V_k \otimes V \cong V_{k+1} \oplus V_{k-1}, \quad k \geq 2.
\end{gather*}
\end{Proposition}

\begin{proof}
Let $V_k$ be the submodule generated by the simple tensors $v_1^{\otimes k}, v_2^{\otimes k}$ in $V^{\otimes k}.$ This is a two dimensional submodule, and by observing the action of $a \in \text{Dic}_{\infty}$ on these simple tensors, we can conclude that $V_i \ncong V_j$ for $i < j.$ Furthermore, since $\theta$ is chosen so that $e^{i \theta}$ is not a root of unity, we obtain an infinite family $\{V_k\}_{k \in \N}$ of mutually nonisomorphic two dimensional simples. 

The Clebsch--Gordan type rule can then be proven easily, noting that it is concerned with a four dimensional space. The rule in turn can then be used to show
\begin{align*}
V^{\otimes 2} & \cong V_2 \oplus \bar{\C} \oplus \C, \\ 
V^{\otimes 3} & \cong V_3 \oplus V^{\oplus 3}, \\
V^{\otimes 4} & \cong V_4 \oplus V_2^{\oplus 4} \oplus \bar{\C}^{\oplus 3} \oplus \C^{\oplus 3}, \\
V^{\otimes 5} & \cong V_5 \oplus V_3^{\oplus 5} \oplus V^{\oplus 10},
\end{align*}
and so on. In fact, the number of summands grows as $2^n$ by the main theorem of 
\cite{CoOsTu-growth} and the fusion graph that the reader needs to know to extract the above rules (the summands of $V^{\otimes k}$ correspond 1:1 to walks of length $k$ in this graph) is the infinite (affine) type D Dynkin diagram (with black affine node):
\begin{gather*}
\tilde{D}_{\infty}=D_{\infty}=
\lim_{n\to\infty}
\dynkin[
labels={0',0,1,2,,n{-}2,n{-}1,n',n},
label directions={,,left,,,,right,,},
scale=1.8,
extended
] D{oooo...oo*o}
.
\end{gather*}
This together with \autoref{lemma:all-appear} shows that 
every finite dimensional simple representation of the dicyclic group is either trivial, $\bar{\C}$, or is two dimensional and isomorphic to $V_k$ for some $k \in \N$. Note also that $V_j$ is not a summand of $V^{\otimes k}$ if $j<k$: that is, each tensor power $V^{\otimes k}$ decomposes into $V_k$ and some summands previously encountered.
\end{proof}

\begin{Remark}
As stated in the proof above, the number of summands grows as $2^n$. This should lead 
to interesting examples in the sense of \cite{KhSiTu-monoidal-cryptography}.
\end{Remark}

\begin{Lemma}\label{L:JWscalar}
For all $k\in\N$, the tensor power $V^{\otimes k}$ contains 
a multiplicity one summand, the \emph{leading summand}, which has not appeared for a smaller tensor power. 
The ratios of the dimensions of these leading summands is
\begin{gather*}
1,2,1,1,1,\dots
.
\end{gather*}
\end{Lemma}

\begin{proof}
Directly from \autoref{P:cgrule}.
\end{proof}

\begin{Lemma}
We have $V_n \otimes \bar{\C} \cong V_n$.
\end{Lemma}

\begin{proof}
Since $a$ acts as the identity on $\bar{\C}$.
\end{proof}

Thus, the decompositions of the tensor powers of the natural representation of the dicyclic group follow a pattern similar to that of $SU(2)$, whose fusion graph for tensoring with $V$ is
\begin{gather*}
A_{\infty}=\lim_{n\to\infty}
\dynkin[labels={1,2,,n{-}1,n}]A5
,
\end{gather*}
but with an additional new one dimensional summand $\bar{\C}$ in even tensor powers.

We now continue with the representation theory of the 
finite dicyclic groups.

\begin{Proposition}\label{P:fusion-graph}
We have the following.
\begin{enumerate}

\item There are precisely $n+3$ nonisomorphic simple $\text{Dic}_n$-representations: 
\begin{gather*}
\text{1 dimensional}\colon \{V_0=\C,V_{0'}=\bar{\C},V_{n},V_{n'}\},\quad
\text{2 dimensional}\colon \big\{V_k|k \in \{1,\dots,n-1\}\big\}.
\end{gather*}

\item $\C$ and $\bar{\C}$ are as before, while, for $k = 1,2, \dots, n-1$, all of the two dimensional simple representations $V_k$ are given by
\begin{gather*}
\begin{NiceMatrixBlock}[auto-columns-width]
x \mapsto \begin{pNiceMatrix}
0 & -1 \\[0.5em]
1 & 0
\end{pNiceMatrix}, \quad
a \mapsto \begin{pNiceMatrix}
e^{\frac{k \pi i}{n}} & 0 \\[0.5em] 
0 & e^{\frac{-k \pi i}{n}}
\end{pNiceMatrix}
\end{NiceMatrixBlock}
.
\end{gather*}
Moreover, for $k=n$ the corresponding representation is $V_{n}\oplus V_{n'}$.

\item Let $V$ be the $\text{Dic}_n$-representation induced from the 
$\text{Dic}_\infty$-representation on $V=\C^2$. Then all the above simple $\text{Dic}_n$-representations appear as direct summands of some $V^{\otimes k}$.

\item The fusion graph of tensoring with $V$ is the affine type D Dynkin diagram:
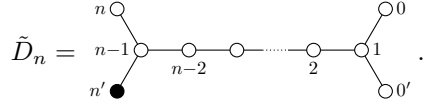
\begin{gather*}
\tilde{D}_{n}=
\dynkin[
labels={0',0,1,2,,n{-}2,n{-}1,n',n},
label directions={,,left,,,,right,,},
scale=1.8,
extended
] D{oooo...oo*o}.
\end{gather*}

\end{enumerate}
\end{Proposition}

\begin{proof}
The degree $n$ dicyclic group $\text{Dic}_n$ consists of $n+3$ conjugacy classes. The center is a union of the singleton classes $\{1\}, \{x^2\}$. There are $n-1$ classes of size two, given by the elements $\{a,a^{-1}\}, \{a^2, a^{-2}\}$, up to $\{a^{n-1}, a^{-n+1}\}$, and finally two classes with $n$ elements: $\{x,xa^2, \dots, xa^{2n-2}\}, \{xa, xa^3, \dots, xa^{2n-1}\}.$ Hence, the count comes out as $n+3$ simple $\text{Dic}_n$-representations. One can then check that the given matrix representations are pairwise nonisomorphic and simple, and the claim follows.

Alternatively, we briefly outline the construction of all the simples as summands of tensor powers of the natural representation, in the same way as above for $\text{Dic}_{\infty}$. For $1 \leq k < n$, the decomposition is completely identical: the two dimensional simple $V_k$ appears as a multiplicity one summand of $V^{\otimes k}$. Indeed, if $k < n$, the simple tensors $\{v_1^{\otimes k}, v_2^{\otimes k}\}$ generate the two dimensional summand $V_k \subset V^{\otimes k}$. There are two new one dimensional simple representations for $\text{Dic}_n$. Observe that for $k = n$, the element $v_1^{\otimes k} + v_2^{\otimes k}$ is acted on by a sign by $a \in G$ and fixed by $x \in G$ when $n$ is even and hence generates a 1-dimensional simple not isomorphic to $\bar{\C}.$ In the case of $n$ odd, one takes $v_1^{\otimes k} - v_2^{\otimes k}$. Denote this simple representation by $V_n$. Finally, also inside of $V^{\otimes k}$ for $k = n$ with $n$ even, the element $v_1^{\otimes k} - v_2^{\otimes k}$ is transformed by a sign under both $a,x \in G$ and hence generates another one dimensional simple, not isomorphic to $V_n, \bar{\C}$, or the trivial representation: we denote it by $V_{n'}.$ Again, if $n$ is odd, one needs to change a sign on the generating element, but $V_{n'}$ remains a summand.
\end{proof}

As outlined in the introduction, we would like to study the endomorphism algebra $A_k = \text{End}_{G}(V^{\otimes k})$ for both the dicyclic group $\text{Dic}_{\infty}$ and its finite quotient group $G = \text{Dic}_n$, and compare it to the Temperley--Lieb algebra $TL_k(-2).$

\begin{Lemma}\label{lemma:dimension-dicinf}
Let $V$ denote the standard representation of $G = \mathrm{Dic}_{\infty}$ and $\chi$ its character. For all integers $k\in\Z_{\geq 1}$, we have
\begin{gather*}
\dim_{\C} \mathrm{End}_{G}(V^{\otimes k}) = \langle 1, \chi^{2k} \rangle = \frac{2^{2k}}{4\pi}\displaystyle\int_{0}^{2\pi} \cos^{2k}(\theta) \ d\theta = \frac{1}{2}\binom{2k}{k}.
\end{gather*}
\end{Lemma}

\begin{proof}
First, we have $\dim_{\C}\text{End}_{G}(V^{\otimes k}) = \dim_{\C} \text{Hom}_G(\C, V^{\otimes 2k})$ since $V$ is self dual. To justify the first equality, note that as a discrete subgroup of $SU(2)$, $G$ acts continuously on $V^{\otimes 2k}$ and hence the left-hand side equals $\dim \text{Hom}_{\bar{G}}(\C, V^{\otimes 2k})$. Since $\bar{G}$ is compact, we can compute the dimension with a character calculation by integrating the character $\chi^{2k}$ over all of $\bar{G}$ and normalizing:
\begin{gather*}
\langle 1, \chi^{2k} \rangle = \frac{2^{2k}}{4\pi}\displaystyle\int_{0}^{2\pi} \cos^{2k}(\theta).
\end{gather*}

Finally, the identity $\frac{2^{2k}}{4\pi}\int_{0}^{2\pi} \cos^{2k}(\theta) \ d\theta = \frac{1}{2}\binom{2k}{k}$ may be proven by converting to an integral over the unit circle in the complex plane and applying the Cauchy residue theorem.
\end{proof}

Let us now take the first steps to identify the dimension of $A_k$ in the case of the degree $n$ dicyclic group $G=\text{Dic}_n$. 

\begin{Proposition}\label{P:diminf}
For all integers $1 \leq k < n$, we have
\begin{gather*}
\dim_{\C} \mathrm{End}_{G}(V^{\otimes k}) = \langle 1, \chi^{2k} \rangle = \frac{2^{2k}}{4\pi}\displaystyle\int_{0}^{2\pi} \cos^{2k}(\theta) \ d\theta = \frac{1}{2}\binom{2k}{k}.
\end{gather*}
\end{Proposition}

\begin{proof}
As $V$ is self-dual, $\dim \text{End}_{G}(V^{\otimes k})$ equals the multiplicity of the trivial representation in $V^{\otimes 2k}$. Now, use
\autoref{lemma:dimension-dicinf} and its proof.
\end{proof}

Note that $\dim_{\C} \mathrm{End}_{G}(V^{\otimes n}) = \frac{1}{2}\binom{2n}{n} + 1$, which is the case when we hit the second fish tail of the affine type D Dynkin diagram. The extra morphism in the $k+1$ endomorphism algebra arises as the primitive idempotent projecting to the summand $V_{n'}$. In the case of $\text{Dic}_{\infty}$, this summand does not exist: every tensor power of $V$ results in a unique two dimensional summand of multiplicity one generated by the vectors $v_1^{\otimes k}, v_2^{\otimes k}$. In the case of the finite dicyclic group $\text{Dic}_n$, the splitting arises because the order of $a \in G$ is finite. In particular, the breaking point $k = n$ comes from the exponential $e^{ \pi i/n}$ in the matrix for the generator $a \in G \subset SU(2).$

One would expect the dimension of $A_k$ to be $\frac{1}{2}\binom{2k}{k}$ for all $k$ in the case of the dicyclic group $\text{Dic}_{\infty}$, since there is no ``breaking point'' in the tensor power decompositions. To compute the dimension of $A_k$ in the case of the finite dicyclic group $\text{Dic}_n$, we prove the following, where we write the character table using the conjugacy classes as index for the columns, with value vectors denoted by $col_{index}$. Let $[2]_{d,n}=q+q^{-1}$ for $q=\exp(di\pi/n)$.

\begin{Theorem}\label{characterformula}
For all integers $k\in\Z_{\geq 1}$, we have
\begin{gather*}
\dim_{\C} \mathrm{End}_{G}(V^{\otimes k})
=
\bigg(\sum\Big(\tfrac{1}{4n}\cdot
\big(1\cdot col_{1}\cdot 2^k
+1\cdot col_{x^2}\cdot (-2)^k
+\sum_{d=1}^{n-1}2\cdot col_{a^d,a^{-d}}\cdot[2]_{d,n}^k\big)
\Big)\bigg)^2,
\end{gather*}
where the outside sum means summing over all entries of the vector. Moreover, the vectors $col_{index}$ can be explicitly obtained from \autoref{L:Dicyclic}.(b) and the character table of the dihedral group as in \autoref{L:DihedralCharacterTable}.
\end{Theorem}

\begin{proof}
Recall that the fusion graph we need to study is the affine type D Dynkin diagram, see 
\autoref{P:fusion-graph}. We need to take the $k$th power of the adjacency matrix of this graph, and compute the sum of the squares of the entries in the first row.
By standard argument, see for example \cite{EtGeNiOs-tensor-categories,LaTuVa-growth-pfdim,LaTuVa-growth-pfdim-inf} and in particular \cite{He}, this boils down to a character calculation, and we need the following lemma. To this end, recall that we explicitly computed the conjugacy classes of the dicyclic group in the proof of \autoref{P:fusion-graph}.

\begin{Lemma}\label{L:DihedralCharacterTable}
We have the following character tables for the dihedral group.
Let $\ran=2\pi/n$.
Firstly, if $n$ is even then:
\begin{center}
\scalebox{0.8}{\begin{tabular}{C||CCCCCCCCC}
& 1 & a & a^{2} & a^{3} & \dots & a^{m-1} & a^{m} & b & ab \\
& 1 & 2 & 2 & 2 & \dots & 2 & 1 & m & m \\
\hline
\hline
\chi_{0} & 1 & 1 & 1 & 1 & \dots & 1 & 1 & 1 & 1 \\
\chi_{0}^{\star} & 1 & 1 & 1 & 1 & \dots & 1 & 1 & -1 & -1 \\
\chi_{1} & 2 & 2\cos(\ran) & 2\cos(2\ran) & 2\cos(3\ran) & \dots & 2\cos\big((m-1)\ran\big) & -2 & 0 & 0 \\
\vdots & \vdots & \vdots & \vdots & \vdots & \dots & \vdots & \vdots & \vdots & \vdots \\
\chi_{m-1} & 2 & 2\cos\big((m-1)\ran\big) & 2\cos\big(2(m-1)\ran\big) & 2\cos\big(3(m-1)\ran\big) & \dots & 2\cos\big((m-1)^{2}\ran\big) & (-1)^{m-1}2 & 0 & 0 \\
\chi_{m} & 1 & -1 & 1 & -1 & \dots & (-1)^{m-1} & (-1)^{m} & -1 & 1 \\
\chi_{m}^{\star} & 1 & -1 & 1 & -1 & \dots & (-1)^{m-1} & (-1)^{m} & 1 & -1 \\
\end{tabular}}
.
\end{center}
Second, if $n$ is odd, then:
\begin{center}
\scalebox{0.8}{\begin{tabular}{C||CCCCCCCC}
& 1 & a & a^{2} & a^{3} & \dots & a^{m-1} & a^{m} & b \\
& 1 & 2 & 2 & 2 & \dots & 2 & 2 & n \\
\hline
\hline
\chi_{0} & 1 & 1 & 1 & 1 & \dots & 1 & 1 & 1 \\
\chi_{0}^{\star} & 1 & 1 & 1 & 1 & \dots & 1 & 1 & -1 \\
\chi_{1} & 2 & 2\cos(\ran) & 2\cos(2\ran) & 2\cos(3\ran) & \dots & 2\cos\big((m-1)\ran\big) & 2\cos(m\ran) & 0 \\
\vdots & \vdots & \vdots & \vdots & \vdots & \dots & \vdots & \vdots & \vdots \\
\chi_{m} & 2 & 2\cos(m\ran) & 2\cos(2m\ran) & 2\cos(3m\ran) & \dots & 2\cos\big((m-1)m\ran\big) & 2\cos(m^{2}\ran) & 0 \\
\end{tabular}}
.
\end{center}
Here $m$ is determined by $n = 2m$ or $n = 2m + 1$. The first row is a representative of the associated conjugacy class, and the second row is the number of elements in the class.
\end{Lemma}

\begin{proof}
This is well-known, e.g. this follows from \cite[Section 5.3]{Se-rep-theory-finite-groups}, and omitted. See also the additional file for the paper \cite{GiTuWi-pl-reps} where this is stolen from.
\end{proof}

\begin{Lemma}\label{L:chartable}
The character values $\chi$ on $V$ are as follows:
\begin{enumerate}

\item The conjugacy class $\{1\}$ takes value $\chi=2$.

\item The conjugacy class $\{x^2\}$ takes value $\chi=-2$.

\item The conjugacy class $\{a^d,a^{-d}\}$ takes value $\chi=[2]_{d,n}$.

\item The conjugacy classes with $n$ elements take values $\chi=0$.

\end{enumerate}
\end{Lemma}

\begin{proof}
The character table of the dihedral groups is well known; see the proof of \autoref{L:DihedralCharacterTable} above.
We combine this with \autoref{L:Dicyclic}.(b).
\end{proof}

Following \cite{He} and using \autoref{L:chartable}, one gets the desired result.
\end{proof}

\begin{Theorem}\label{T:CharFormula}
We have
\begin{gather*}
\dim_{\C} \mathrm{End}_{G}(V^{\otimes k}) = \frac{2^{2k+1}}{4n}\left( 1 + \sum\limits_{m=1}^{n-1} \cos^{2k}\left(\frac{m\pi}{n}\right)\right).
\end{gather*}
\end{Theorem}

\begin{proof}
Carrying out the above character calculation explicitly, the formula of \autoref{characterformula} gives us the claimed result.
\end{proof}

\section{Diagrammatics for dicyclic groups}\label{S:Dia}

We now describe the representation category of the dicyclic group diagrammatically.
Our to-go reference for Temperley--Lieb diagrammatics is \cite{KaLi-TL-recoupling}, where
the reader can find many facts that we implicitly use.

\subsection{The dicyclic Temperley--Lieb algebra $DiTL_k(\delta)$}

To proceed towards our goal of modifying the Temperley--Lieb category for $SU(2)$, we must now more concretely identify the structure of the endomorphism algebras $A_k$.

We start this section with the construction of what we are calling the \emph{dicyclic Temperley--Lieb algebras}, which are quotients of the type D Temperley--Lieb algebras. We will recall the definition and diagrammatic interpretation of these algebras here, but for a more thorough treatment see \cite{Gr-dots-tl} or \cite{LiSh-tl-typed}. The type D Dynkin diagram, which we now need instead of the affine one, is:
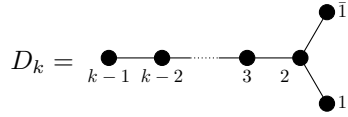
\begin{gather*}
D_{k}=
\reflectbox{\dynkin[scale=2, labels={\reflectbox{$k-1$}, \reflectbox{$k-2$}, \reflectbox{$3$}, \reflectbox{$2$}, \reflectbox{$\bar{1}$}, \reflectbox{$1$}}] D{oo...oooo}}
.
\end{gather*}
For each vertex labeled $1,\bar{1}, 2, \dots, k-1$ in the diagram, we associate a generator $e_i$. 

\begin{Definition}
The \emph{type D Temperley--Lieb algebra} is the free $\C$-algebra generated by (generators indexed by the vertices of $D_k$) $e_1, e_{\bar{1}}, e_2, \dots, e_{k-1}$, subject to the relations:
\begin{enumerate}

\item $e_i^2 = \delta e_i$ for all $i \in \{1, \bar{1}, 2, \dots, k-1\}.$

\item $e_ie_j = e_je_i$ if $i,j$ are disconnected in the diagram.

\item $e_ie_je_i = e_i$ if $i,j$ are connected in the diagram.

\end{enumerate}
Here, as throughout, we fix the so-called circle parameter $\delta\in\C$.
\end{Definition}

\begin{Remark}
We could also work over $\C[\delta]$ and make $\delta$ a parameter. Since we will only need $\delta=-2$ the above formulation is sufficient for us.
\end{Remark}

This is the usual construction of the Temperley--Lieb algebra for the type A diagram, but with a single additional node $\bar{1}$:

\begin{Lemma}\label{L:oldTL}
The subalgebra generated by $e_1, e_2, \dots, e_{k-1}$ is the usual Temperley--Lieb algebra $TL_k(\delta).$
\end{Lemma}

\begin{proof}
Not difficult (and not needed); here is a sketch.
Let $A$ be the subalgebra generated by $e_1,e_2,\dots,e_{k-1}$.
Restricting the defining relations (a)--(c) to the vertices $1,2,\dots,k-1$ (which form the type $A_{k-1}$ string) gives precisely the usual Temperley--Lieb relations.
Hence, the assignment $U_i\mapsto e_i$ induces a surjective algebra homomorphism
$TL_k(\delta)\twoheadrightarrow A$.

For injectivity, we can use the results in \autoref{S:Inf} to see that both algebras admit a representation on $\mathrm{End}_{\C}(V^{\otimes k})$ by the same matrices.
\end{proof}

Diagrammatically, the generators $e_i$ may be thought of as the usual cup-cap generators of the Temperley--Lieb algebra, cf. 
\autoref{L:oldTL}, together with the following \emph{decorated diagram} for the generator $e_{\bar{1}}$ corresponding to the new node:
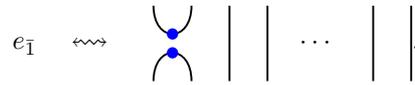
\begin{gather*}
e_{\bar{1}} \quad
\leftrightsquigarrow
\begin{tikzpicture}[anchorbase, scale=0.5]
\draw [draw=none] (0,0) rectangle (3,2);
\draw [usual] [postaction={decorate,decoration={markings,mark=at position 0.5 with {\node[circle,fill=blue,inner sep=1.5pt] {};}}}] (1,0) [dndn=1];
\draw [usual] [postaction={decorate,decoration={markings,mark=at position 0.5 with {\node[circle,fill=blue,inner sep=1.5pt] {};}}}] (1,2) [upup=1];
\draw [usual] (3,0) [dnup=0];
\draw [usual] (4,0) [dnup=0];
\end{tikzpicture} \ \ 
\ \cdots
\begin{tikzpicture}[anchorbase, scale=0.5]
\draw [draw=none] (0,0) rectangle (2,2);
\draw [usual] (1,0) [dnup=0];
\draw [usual] (2,0) [dnup=0];
\end{tikzpicture}
.
\end{gather*}
When interpreted this way, we have the following rules for \emph{decoration removal} when we multiply diagrams:
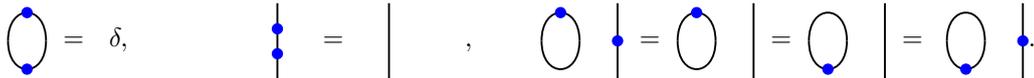
\vspace{-10pt}
\begin{gather*}
\begin{tikzpicture}[anchorbase, scale=0.5]
\draw [draw=none] (0,0) rectangle (2,4);
\draw [usual] [postaction={decorate,decoration={markings,mark=at position 0.5 with {\node[circle,fill=blue,inner sep=1.5pt] {};}}}] (1,2) [upup=1];
\draw [usual] [postaction={decorate,decoration={markings,mark=at position 0.5 with {\node[circle,fill=blue,inner sep=1.5pt] {};}}}] (1,2) [dndn=1];
\end{tikzpicture} \ = \  \delta, \quad 
\begin{tikzpicture}[anchorbase, scale=0.5]
\draw [draw=none] (0,0) rectangle (2,2);
\draw [usual] [postaction={decorate,decoration={markings,mark=at position 0.66 with {\node[circle,fill=blue,inner sep=1.5pt] {};}}}] [postaction={decorate,decoration={markings,mark=at position 0.33 with {\node[circle,fill=blue,inner sep=1.5pt] {};}}}] (1,0) [dnup=0];
\end{tikzpicture}
=
\begin{tikzpicture}[anchorbase, scale=0.5]
\draw [draw=none] (0,0) rectangle (3,2);
\draw [usual] (1,0) [dnup=0];
\end{tikzpicture}\hspace{-0.95cm}, \quad
\begin{tikzpicture}[anchorbase, scale=0.5]
\draw [draw=none] (0,0) rectangle (2,4);
\draw [usual] (1,2) [upup=1];
\draw [usual] [postaction={decorate,decoration={markings,mark=at position 0.5 with {\node[circle,fill=blue,inner sep=1.5pt] {};}}}] (1,2) [dndn=1];
\draw [usual] [postaction={decorate,decoration={markings,mark=at position 0.5 with {\node[circle,fill=blue,inner sep=1.5pt] {};}}}] (3,1) [dnup=0];
\end{tikzpicture} \ = \
\begin{tikzpicture}[anchorbase, scale=0.5]
\draw [draw=none] (1,0) rectangle (2,4);
\draw [usual] (1,2) [upup=1];
\draw [usual] [postaction={decorate,decoration={markings,mark=at position 0.5 with {\node[circle,fill=blue,inner sep=1.5pt] {};}}}] (1,2) [dndn=1];
\draw [usual] (3,1) [dnup=0];
\end{tikzpicture}
\ = \
\begin{tikzpicture}[anchorbase, scale=0.5]
\draw [draw=none] (1,0) rectangle (2,4);
\draw [usual][postaction={decorate,decoration={markings,mark=at position 0.5 with {\node[circle,fill=blue,inner sep=1.5pt] {};}}}] (1,2) [upup=1];
\draw [usual] (1,2) [dndn=1];
\draw [usual] (3,1) [dnup=0];
\end{tikzpicture}
\ = \ \hspace{-0.4cm}
\begin{tikzpicture}[anchorbase, scale=0.5]
\draw [draw=none] (0,0) rectangle (2,4);
\draw [usual][postaction={decorate,decoration={markings,mark=at position 0.5 with {\node[circle,fill=blue,inner sep=1.5pt] {};}}}] (1,2) [upup=1];
\draw [usual] (1,2) [dndn=1];
\draw [usual] [postaction={decorate,decoration={markings,mark=at position 0.5 with {\node[circle,fill=blue,inner sep=1.5pt] {};}}}] (3,1) [dnup=0];
\end{tikzpicture}
.
\end{gather*}
Note that in the type $D$ Temperley--Lieb algebra, loops with one decoration do not evaluate to $\delta.$ This means that a basis will include some diagrams with loops. However, if we further impose the relation $e_1e_{\bar{1}} = e_{\bar{1}}e_1 = 0$ on the type D Temperley--Lieb algebra, the third decoration removal relation vanishes, and we get:
\begin{gather*}
\begin{tikzpicture}[anchorbase, scale=0.5]
\draw [draw=none] (1,0) rectangle (2,4);
\draw [usual] (1,2) [upup=1];
\draw [usual] [postaction={decorate,decoration={markings,mark=at position 0.5 with {\node[circle,fill=blue,inner sep=1.5pt] {};}}}] (1,2) [dndn=1];
\end{tikzpicture}
=
\begin{tikzpicture}[anchorbase, scale=0.5,yscale=-1]
\draw [draw=none] (1,0) rectangle (2,4);
\draw [usual] (1,2) [upup=1];
\draw [usual] [postaction={decorate,decoration={markings,mark=at position 0.5 with {\node[circle,fill=blue,inner sep=1.5pt] {};}}}] (1,2) [dndn=1];
\end{tikzpicture}
=0.
\end{gather*}

\begin{Example}
The diagram corresponding to the element $e_{\bar{1}}e_2e_1$ in the 3 strand algebra is: 
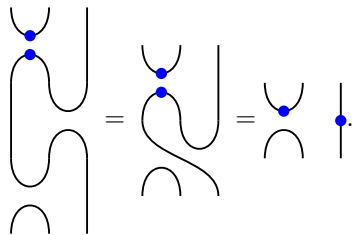
\begin{gather*}
\begin{tikzpicture}[anchorbase, scale=0.5]
\draw [draw=none] (1,0) rectangle (2,4);
\draw [usual] (2,2) [upup=1];
\draw [usual] (1,0) [dnup=0];
\draw [usual] [postaction={decorate,decoration={markings,mark=at position 0.5 with {\node[circle,fill=blue,inner sep=1.5pt] {};}}}] (1,4) [upup=1];
\draw [usual] [postaction={decorate,decoration={markings,mark=at position 0.5 with {\node[circle,fill=blue,inner sep=1.5pt] {};}}}] (1,2) [dndn=1];
\draw [usual] (3,2) [dnup = 0];
\draw [usual] (2,0) [dndn=1];
\draw [usual] (3,-2) [dnup=0];
\draw [usual] (1,-2) [dndn=1];
\draw [usual] (1,0) [upup=1];
\end{tikzpicture} \ = \
\begin{tikzpicture}[anchorbase, scale=0.5]
\draw [draw=none] (1,0) rectangle (2,4);
\draw [usual] [postaction={decorate,decoration={markings,mark=at position 0.5 with {\node[circle,fill=blue,inner sep=1.5pt] {};}}}] (1,4) [upup=1];
\draw [usual] [postaction={decorate,decoration={markings,mark=at position 0.5 with {\node[circle,fill=blue,inner sep=1.5pt] {};}}}] (1,2) [dndn=1];
\draw [usual] (3,2) [dnup = 0];
\draw [usual] (2,2) [upup = 1];
\draw [usual] (3,0) [dnup = -2];
\draw [usual] (1,0) [dndn = 1];
\end{tikzpicture}   \ = 
\begin{tikzpicture}[anchorbase, scale=0.5]
\draw [draw=none] (1,0) rectangle (2,2);
\draw [usual] (1,0) [dndn=1];
\draw [usual] [postaction={decorate,decoration={markings,mark=at position 0.5 with {\node[circle,fill=blue,inner sep=1.5pt] {};}}}] (1,2) [upup=1];
\draw [usual] [postaction={decorate,decoration={markings,mark=at position 0.5 with {\node[circle,fill=blue,inner sep=1.5pt] {};}}}] (3,0) [usual] [dnup=0];
\end{tikzpicture},
\end{gather*}
where we used the standard zigzag relations of the Temperley--Lieb algebra.
\end{Example}

We now define one of our main players:

\begin{Definition}
The \emph{dicyclic Temperley--Lieb algebra} on $k$ strands, denoted $DiTL_k(\delta)$, is the quotient of the type D Temperley--Lieb algebra by the two-sided ideal generated by $e_1e_{\bar{1}}$.
\end{Definition}

\begin{Lemma}\label{L:tlbasis}
The dimension of the dicyclic Temperley--Lieb algebra $DiTL_k(\delta)$ is $\frac{1}{2}\binom{2k}{k}$. The algebra has a basis consisting of diagrams with satisfying the following two conditions:
\begin{enumerate}

\item There are an even number of decorations.

\item Any decorated strand is exposed to the left boundary of the diagram.

\end{enumerate}
\end{Lemma}

\begin{proof}
This is \cite[Lemma 6.5]{Gr-dots-tl}.
\end{proof}

The second condition in \autoref{L:tlbasis}, in other words, says that any decorated strand can be pulled out from the left-hand side of the diagram.

\begin{Example}
For example, the following ten diagrams constitute a basis for $DiTL_3(\delta):$
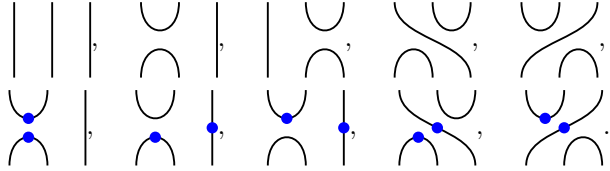
\begin{gather*}
\begin{tikzpicture}[anchorbase,scale=0.5]
\draw [draw=none] (0,0) rectangle (3,2);
\draw [usual] (1,0) [dnup=0];
\draw [usual] (2,0) [dnup=0];
\draw [usual] (3,0) [dnup=0];
\end{tikzpicture},
\begin{tikzpicture}[anchorbase,scale=0.5]
\draw [draw=none] (0,0) rectangle (3,2);
\draw [usual] (1,0) [dndn=1];
\draw [usual] (1,2) [upup=1];
\draw [usual] (3,0) [dnup=0];
\end{tikzpicture},
\begin{tikzpicture}[anchorbase,scale=0.5]
\draw [draw=none] (0,0) rectangle (3,2);
\draw [usual] (1,0) [dnup=0];
\draw [usual] (2,2) [upup=1];
\draw [usual] (2,0) [dndn=1];
\end{tikzpicture},
\begin{tikzpicture}[anchorbase,scale=0.5]
\draw [draw=none] (0,0) rectangle (3,2);
\draw [usual] (1,0) [dndn=1];
\draw [usual] (3,0) [dnup=-2];
\draw [usual] (2,2) [upup=1];
\end{tikzpicture},
\begin{tikzpicture}[anchorbase,scale=0.5]
\draw [draw=none] (0,0) rectangle (3,2);
\draw [usual] (1,0) [dnup=2];
\draw [usual] (1,2) [upup=1];
\draw [usual] (2,0) [dndn=1];
\end{tikzpicture},
\\
\begin{tikzpicture}[anchorbase,scale=0.5]
\draw [draw=none] (0,0) rectangle (3,2);
\draw [usual] [postaction={decorate,decoration={markings,mark=at position 0.5 with {\node[circle,fill=blue,inner sep=1.5pt] {};}}}] (1,0) [dndn=1];
\draw [usual] [postaction={decorate,decoration={markings,mark=at position 0.5 with {\node[circle,fill=blue,inner sep=1.5pt] {};}}}] (1,2) [upup=1];
\draw [usual] (3,0) [dnup=0];
\end{tikzpicture},
\begin{tikzpicture}[anchorbase,scale=0.5]
\draw [draw=none] (0,0) rectangle (3,2);
\draw [usual] [postaction={decorate,decoration={markings,mark=at position 0.5 with {\node[circle,fill=blue,inner sep=1.5pt] {};}}}] (1,0) [dndn=1];
\draw [usual] (1,2) [upup=1];
\draw [usual] [postaction={decorate,decoration={markings,mark=at position 0.5 with {\node[circle,fill=blue,inner sep=1.5pt] {};}}}] (3,0) [dnup=0];
\end{tikzpicture},
\begin{tikzpicture}[anchorbase,scale=0.5]
\draw [draw=none] (0,0) rectangle (3,2);
\draw [usual] (1,0) [dndn=1];
\draw [usual] [postaction={decorate,decoration={markings,mark=at position 0.5 with {\node[circle,fill=blue,inner sep=1.5pt] {};}}}] (1,2) [upup=1];
\draw [usual] [postaction={decorate,decoration={markings,mark=at position 0.5 with {\node[circle,fill=blue,inner sep=1.5pt] {};}}}] (3,0) [usual] [dnup=0];
\end{tikzpicture},
\begin{tikzpicture}[anchorbase,scale=0.5]
\draw [draw=none] (0,0) rectangle (3,2);
\draw [usual] [postaction={decorate,decoration={markings,mark=at position 0.5 with {\node[circle,fill=blue,inner sep=1.5pt] {};}}}] (1,0) [dndn=1];
\draw [usual] [postaction={decorate,decoration={markings,mark=at position 0.5 with {\node[circle,fill=blue,inner sep=1.5pt] {};}}}] (3,0) [usual] [dnup=-2];
\draw [usual] (2,2) [upup=1];
\end{tikzpicture},
\begin{tikzpicture}[anchorbase,scale=0.5]
\draw [draw=none] (0,0) rectangle (3,2);
\draw [usual] [postaction={decorate,decoration={markings,mark=at position 0.5 with {\node[circle,fill=blue,inner sep=1.5pt] {};}}}] (1,0) [dnup=2];
\draw [usual] [postaction={decorate,decoration={markings,mark=at position 0.5 with {\node[circle,fill=blue,inner sep=1.5pt] {};}}}] (1,2) [upup=1];
\draw [usual] (2,0) [dndn=1];
\end{tikzpicture}
.
\end{gather*}
Note that there are only zero or two dots.
\end{Example}

\begin{Remark}
We note that these algebras are also called the \emph{forked} Temperley--Lieb algebras due to the fish tail of the Dynkin diagram. 
\end{Remark}

\subsection{The dicyclic Temperley--Lieb category $\mathbf{DiTL}(\delta, \infty)$}

The reader might now want to recall some basics of monoidal categories, see e.g. \cite{EtGeNiOs-tensor-categories,Tu-qt}.

\begin{Definition}
The \emph{Temperley--Lieb category} $\mathbf{TL}(\delta)$, is the free $\C$-linear pivotal category generated by a single object $\bullet$ of 
categorical dimension $-2$. The relations thus are:
\begin{equation*}
(1) \ \text{ev} \circ \text{coev} = \delta \cdot \text{id}, \quad \quad \quad (2) \ (\text{ev} \otimes \text{id}) \circ (\text{id} \otimes \text{coev}) = (\text{coev} \otimes \text{id}) \circ (\text{id} \otimes \text{ev}) = \text{id},
\end{equation*}
The morphisms $\text{ev}$ and $\text{coev}$ are the ones coming from the pivotal structure of $\mathbf{TL}(\delta)$.
\end{Definition}

\noindent To properly define a dicyclic version of the Temperley\textendash Lieb category, recall from \autoref{L:tlbasis} that deocrations must always reachable from the left, e.g.:
\begin{gather*}
\text{allowed}:
\begin{tikzpicture}[anchorbase,scale=0.5]
\draw [draw=none] (0,0) rectangle (3,2);
\draw [usual][postaction={decorate,decoration={markings,mark=at position 0.5 with {\node[circle,fill=blue,inner sep=1.5pt] {};}}}] (1,0) [dnup=0];
\draw [usual] (2,0) [dnup=0];
\end{tikzpicture}
,\quad
\text{not allowed}:
\begin{tikzpicture}[anchorbase,scale=0.5]
\draw [draw=none] (0,0) rectangle (3,2);
\draw [usual] (1,0) [dnup=0];
\draw [usual][postaction={decorate,decoration={markings,mark=at position 0.5 with {\node[circle,fill=blue,inner sep=1.5pt] {};}}}] (2,0) [dnup=0];
\end{tikzpicture}
.
\end{gather*}
This suggests that we need to work with (right) module categories. 
Recall a right module category over a monoidal category, see e.g. 
\cite[Section 2]{HaOl-actions-tensor-categories} or many other references 
such as \cite{We-intro-hom-alg,EtGeNiOs-tensor-categories,SaTu-bcd-webs}.
The categorical version of the above is:

\begin{Definition}
The \emph{dicyclic Temperley--Lieb category} $\mathbf{DiTL}(\delta, \infty)$, is the free $\C$-linear $\mathbf{TL}(\delta)$-module category generated by a morphism $e_{\bar{1}} \in \text{End}(\bullet^{\otimes 2})$ subject to the following relations:
\begin{gather*}
(3) \ (\text{id} \otimes \text{ev}) \circ (e_{\bar{1}} \otimes \text{id}) \circ (\text{id} \otimes \text{coev}) = \text{id}, \quad \quad \quad (4) \ (\text{ev} \otimes \text{id}) \circ (\text{id} \otimes e_{\bar{1}}) \circ (\text{coev} \otimes \text{id}) = \text{id},
\\
(5) \ e_{\bar{1}} \circ \text{ev} = \text{coev} \circ e_{\bar{1}} = 0, \quad \quad \quad (6) \ e_{\bar{1}}^2 = \delta \cdot e_{\bar{1}}.
\end{gather*}
Note that we also have (1) and (2) from $\mathbf{TL}(\delta)$, by definition.
\end{Definition}

The relations (1) and (2) are those of the Temperley--Lieb category (where we preferred the notation using evaluations and coevaluations over $e_i$), while relations (3)-(6) encode the properties of the generator $e_{\bar{1}}$ in $DiTL_n(\delta).$ Loops with two decorations evaluate to a scalar just like in the Temperley--Lieb category, however, by equation (5), loops with a single decoration vanish.

\begin{Example}
Drawing a decorated cup-cap for the generating morphism $e_{\bar{1}}$, we may interpret the above relations diagrammatically as follows:
\begin{gather*}
\begin{tikzpicture}[anchorbase, scale=0.5]
\draw [draw=none] (1,0) rectangle (2,6);
\draw [usual](2,2) [upup=1];
\draw [usual] (1,0) [dnup=0];
\draw [usual] [postaction={decorate,decoration={markings,mark=at position 0.5 with {\node[circle,fill=blue,inner sep=1.5pt] {};}}}] (1,4) [upup=1];
\draw [usual] [postaction={decorate,decoration={markings,mark=at position 0.5 with {\node[circle,fill=blue,inner sep=1.5pt] {};}}}] (1,2) [dndn=1];
\draw [usual] (2,4) [dndn = 1];
\draw [usual] (3,2) [dnup = 0];
\draw [usual] (1,4) [dnup=0];
\end{tikzpicture}=
\begin{tikzpicture}[anchorbase, scale=0.5]
\draw [draw=none] (1,0) rectangle (1,2);
\draw [usual] [postaction={decorate,decoration={markings,mark=at position 0.33 with {\node[circle,fill=blue,inner sep=1.5pt] {};}}}] [postaction={decorate,decoration={markings,mark=at position 0.66 with {\node[circle,fill=blue,inner sep=1.5pt] {};}}}] (1,0) [dnup = 0];
\end{tikzpicture}=
\begin{tikzpicture}[anchorbase, scale=0.5]
\draw [draw=none] (1,0) rectangle (1,2);
\draw [usual] (1,0) [dnup = 0];
\end{tikzpicture}, \quad \quad \quad \quad 
\quad \begin{tikzpicture}[anchorbase, scale=0.5]
\draw [draw=none] (1,0) rectangle (2,4);
\draw [usual] (1,2) [upup = 1];
\draw [usual] [postaction={decorate,decoration={markings,mark=at position 0.5 with {\node[circle,fill=blue,inner sep=1.5pt] {};}}}] (1,2) [dndn = 1];
\end{tikzpicture} = 0, \quad \quad \quad \quad 
\quad \begin{tikzpicture}[anchorbase, scale=0.5]
\draw [draw=none] (1,0) rectangle (2,4);
\draw [usual] [postaction={decorate,decoration={markings,mark=at position 0.5 with {\node[circle,fill=blue,inner sep=1.5pt] {};}}}] (1,0) [dndn = 1];
\draw [usual] [postaction={decorate,decoration={markings,mark=at position 0.5 with {\node[circle,fill=blue,inner sep=1.5pt] {};}}}] (1,2) [upup = 1];
\draw [usual] [postaction={decorate,decoration={markings,mark=at position 0.5 with {\node[circle,fill=blue,inner sep=1.5pt] {};}}}] (1,2) [dndn = 1];
\draw [usual] [postaction={decorate,decoration={markings,mark=at position 0.5 with {\node[circle,fill=blue,inner sep=1.5pt] {};}}}] (1,4) [upup = 1];
\end{tikzpicture}
=
\ \ \delta \begin{tikzpicture}[anchorbase, scale=0.5]
\draw [draw=none] (0,0) rectangle (3,2);
\draw [usual] [postaction={decorate,decoration={markings,mark=at position 0.5 with {\node[circle,fill=blue,inner sep=1.5pt] {};}}}] (1,0) [dndn=1];
\draw [usual] [postaction={decorate,decoration={markings,mark=at position 0.5 with {\node[circle,fill=blue,inner sep=1.5pt] {};}}}] (1,2) [upup=1];
\end{tikzpicture}.
\end{gather*}
Moreover, relations (1) and (2) take the usual form as the diagrammatic Temperley--Lieb relations.
\end{Example}

\begin{Remark}
Note that the decoration in our diagrams does not represent a generating object of the category $\mathbf{DiTL(\delta,\infty)}$, but we remark here that one could alternatively present this category with an additional generating self dual object of dimension 1 which satisfies the above relations analogously.
\end{Remark}

Define a \emph{crossing map} as follows (we only give the diagrammatic version):
\begin{gather}\label{Eq:crossing}
\begin{tikzpicture}[anchorbase,scale=0.5]
\draw [draw=none] (0,0) rectangle (3,2);
\draw [usual] (0,0) to (2,2);
\draw [usual] (2,0) to (0,2);
\end{tikzpicture}
=
\begin{tikzpicture}[anchorbase,scale=0.5]
\draw [draw=none] (0,0) rectangle (3,2);
\draw [usual] (1,0) [dnup=0];
\draw [usual] (2,0) [dnup=0];
\end{tikzpicture} + 
\begin{tikzpicture}[anchorbase,scale=0.5]
\draw [draw=none] (0,0) rectangle (3,2);
\draw [usual] (1,0) [dndn=1];
\draw [usual] (1,2) [upup=1];
\end{tikzpicture}
.
\end{gather}

We can use this to interpret dots on strands that are not reachable from the left by:
\begin{gather*}
\begin{tikzpicture}[anchorbase,scale=0.5]
\draw [usual] (0,0) to (0,2);
\draw [usual][postaction={decorate,decoration={markings,mark=at position 0.5 with {\node[circle,fill=blue,inner sep=1.5pt] {};}}}] (1,0) to (1,2);
\end{tikzpicture}
=
\begin{tikzpicture}[anchorbase,scale=0.5]
\draw [usual] (0,0) to (0,2);
\draw [usual][postaction={decorate,decoration={markings,mark=at position 0.5 with {\node[circle,fill=blue,inner sep=1.5pt] {};}}}] (1,0) to[out=90,in=270] (-1,1) to[out=90,in=270] (1,2);
\end{tikzpicture}
,
\text{etc.}
\end{gather*}

We can then define a monoidal product on morphisms in $\mathbf{DiTL(\delta,\infty)}$ by juxtaposition.

\begin{Lemma}\label{L:monoidal}
This gives $\mathbf{DiTL(\delta,\infty)}$ the structure of a pivotal category.
\end{Lemma}

\begin{proof}
This follows by interpreting dots as winding around an annulus:
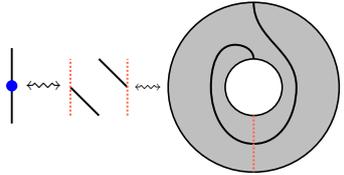
\begin{gather*}
\begin{tikzpicture}[anchorbase,scale=0.5]
\draw [usual][postaction={decorate,decoration={markings,mark=at position 0.5 with {\node[circle,fill=blue,inner sep=1.5pt] {};}}}] (1,0) to (1,2);
\end{tikzpicture}
\leftrightsquigarrow
\scalebox{0.75}{$\begin{tikzpicture}[anchorbase,scale=1]
\draw[usual] (0,0) to (-0.5,0.5);
\draw[usual] (0.5,0.5) to (0,1);
\draw[affine] (-0.5,0) to (-0.5,1);
\draw[affine] (0.5,0) to (0.5,1);
\end{tikzpicture}
\leftrightsquigarrow
\begin{tikzpicture}[anchorbase,scale=1]
\draw[thick,fill=gray!50] (0,0) circle (1.5cm);
\draw[thick,fill=white] (0,0) circle (0.5cm);
\draw[usual] (0,0.5) to[out=90,in=0] (-0.325,0.75) to[out=180,in=90] (-0.75,0) to[out=270,in=180] (0,-1) to[out=0,in=270] (0.75,0)to[out=90,in=270] (0,1.5);
\draw[affine] (0,-0.5) to (0,-1.5);
\end{tikzpicture}$}
\,.
\end{gather*}
In this interpretation the monoidal product becomes
\begin{gather*}
\begin{tikzpicture}[anchorbase,scale=1,yscale=1]
\draw[spinach!35,fill=spinach!35] (-0.5,-0.25) to (-0.5,0.25) to (1.5,0.25) to (1.5,-0.25) to (-0.5,-0.25);
\draw[usual] (0,-1) to (0,-0.25);
\draw[usual] (0,0.25) to (0,1);
\draw[usual] (1,-1) to (1,-0.25);
\draw[usual] (1,0.25) to (1,1);
\node at (0.5,-0.625) {$\dots$};
\node at (0.5,0.625) {$\dots$};
\draw[affine] (-0.5,-1) to (-0.5,1);
\draw[affine] (1.5,-1) to (1.5,1);
\node at (0.5,0) {$f$};
\end{tikzpicture}	
\otimes
\begin{tikzpicture}[anchorbase,scale=1,yscale=1]
\draw[spinach!35,fill=spinach!35] (-0.5,-0.25) to (-0.5,0.25) to (1.5,0.25) to (1.5,-0.25) to (-0.5,-0.25);
\draw[usual] (0,-1) to (0,-0.25);
\draw[usual] (0,0.25) to (0,1);
\draw[usual] (1,-1) to (1,-0.25);
\draw[usual] (1,0.25) to (1,1);
\node at (0.5,-0.625) {$\dots$};
\node at (0.5,0.625) {$\dots$};
\draw[affine] (-0.5,-1) to (-0.5,1);
\draw[affine] (1.5,-1) to (1.5,1);
\node at (0.5,0) {$g$};
\end{tikzpicture}
=
\begin{tikzpicture}[anchorbase,scale=1,yscale=1]
\draw[spinach!35,fill=spinach!35] (-0.5,0.25) to (3.5,0.25) to (3.5,0.75) to (-0.5,0.75) to (-0.5,0.25);
\draw[spinach!35,fill=spinach!35] (-0.5,-0.25) to (3.5,-0.25) to (3.5,-0.75) to (-0.5,-0.75) to (-0.5,-0.25);
\draw[usual] (0,-1) to (0,-0.75);
\draw[usual,crossline] (0,-0.25) to (0,1);
\draw[usual] (1,-1) to (1,-0.75);
\draw[usual,crossline] (1,-0.25) to (1,1);
\draw[usual] (2,-1) to (2,-0.85);
\draw[usual] (2,-0.15) to (2,0.25);
\draw[usual] (2,0.75) to (2,1);
\draw[usual] (3,-1) to (3,-0.85);
\draw[usual] (3,-0.15) to (3,0.25);
\draw[usual] (3,0.75) to (3,1);
\node at (0.5,-1) {$\dots$};
\node at (0.5,1) {$\dots$};
\node at (2.5,-1) {$\dots$};
\node at (2.5,1) {$\dots$};
\draw[affine] (-0.5,-1) to (-0.5,1);
\draw[affine] (3.5,-1) to (3.5,1);
\node at (0.5,-0.5) {$f$};
\node at (2.5,0.5) {$g$};
\end{tikzpicture}
.
\end{gather*}
With this, the statement itself is then a known result, see e.g. \cite{HaOl-actions-tensor-categories,MoSa-affinization,LaTuVa-annular-webs-levi} for details.
\end{proof}

We will from now on use the monoidal and pivotal structure given by \autoref{L:monoidal}.

\section{Infinite dicyclic categories}\label{S:Inf}

We now study the infinite case carefully.

\subsection{An equivalence $\text{Kar}(\mathbf{DiTL}(-2,\infty))  \cong   \mathbf{Rep}(\text{Dic}_{\infty})$}

As before, let $V$ denote the standard representation of the dicyclic group $G = \text{Dic}_{\infty}$, viewed as a subgroup of $SU(2).$ Recall the idempotent $e_{\bar{1}} \in \text{End}_{G}(V^{\otimes 2})$ projecting to the nontrivial one dimensional summand $\bar{\C}$. Let $\mathbf{Rep}(\text{Dic}_{\infty})$ denote the category of finite dimensional complex representations of $\text{Dic}_{\infty}$, and let $\mathcal{C}$ denote the full monoidal subcategory of $\mathbf{Rep}(\text{Dic}_{\infty})$ monoidally generated by $V$.

\begin{Lemma}\label{L:generated}
The morphisms $\text{coev}$, $\text{ev}$, and $e_{\bar{1}}$ monoidally generate $\mathcal{C}$.
\end{Lemma}

\begin{proof}
This follows from \autoref{P:cgrule} and the same arguments as in the Temperley--Lieb case, which can be copied from \cite{RuTeWe-sl2}.
\end{proof}

\begin{Lemma}\label{L:ribbon}
The three categories $\mathbf{DiTL}(-2, \infty)$, $\mathcal{C}$ and $\mathbf{Rep}(\text{Dic}_{\infty})$, are ribbon (with structures specified in the proof).
\end{Lemma}

\begin{proof}
For $\mathbf{DiTL}(-2, \infty)$ we use \autoref{Eq:crossing} as the braiding, and we use the morphisms $\text{ev}$ and $\text{coev}$ for the pivotal structure. For $\mathcal{C}$ and $\mathbf{Rep}(\text{Dic}_{\infty})$ we take the flip map and the duality coming from the Hopf algebra structure as our structures. One easily checks that these are ribbon structures. For $\mathbf{DiTL}(-2, \infty)$ one, for example, calculates:
\begin{gather*}
\begin{tikzpicture}[anchorbase,scale=0.5]
\draw [draw=none] (0,0) rectangle (3,2);
\draw [usual] (0,0) to (2,2);
\draw [usual] (2,0) to (0,2);
\draw [draw=none] (0,2) rectangle (3,4);
\draw [usual] (0,2) to (2,4);
\draw [usual] (2,2) to (0,4);
\end{tikzpicture}
=
\begin{tikzpicture}[anchorbase,scale=0.5]
\draw [draw=none] (0,0) rectangle (3,2);
\draw [usual] (1,0) [dnup=0];
\draw [usual] (2,0) [dnup=0];
\end{tikzpicture}
\text{ since it restricts from $SU(2)$.}
\end{gather*}
The Reidemeister 3 move and the naturality relations can be checked by hand as well, or with {\Magma} if the reader prefers to do that.
\end{proof}

From now on we will use exclusively the ribbon structures as defined in the proof of \autoref{L:ribbon} for the categories in that lemma.

\begin{Lemma}\label{L:functor}
There is a well-defined ribbon $\C$-linear functor $\mathcal{F}\colon\mathbf{DiTL}(-2, \infty) \to\mathcal{C}$ determined by mapping $\bullet \mapsto V$ and the cap, cup, and decorated cup-cap morphisms to $\text{ev}, \text{coev}, $ and $e_{\bar{1}}$ respectively.
\end{Lemma}

\begin{proof}
\autoref{L:relations} shows that $ \text{coev}, \text{ev}$, and $e_{\bar{1}}$ satisfy the defining relations of $\mathbf{DiTL}(-2, \infty)$. The rest of the statement is easy to verify.
\end{proof}

\begin{Theorem}\label{T:equivalence-infinity}
The functor $\mathcal{F}$ from \autoref{L:functor} is fully faithful.
\end{Theorem}

\begin{proof}
In a pivotal category, where ${}^\star$ denotes the dual, we can focus on the endomorphism algebras by using the following, well-known result that we repeat for completeness (the argument works more general, but we will not need this):

\begin{Lemma}\label{L:Bending}
Let $X,Y,Z$ be objects in a pivotal category $C$. Then we have:
\begin{gather*}
\mathrm{Hom}_{C}(X\otimes Y,Z)
\cong
\mathrm{Hom}_{C}\big(Y,
X^\star\otimes Z\big),
\quad
\mathrm{Hom}_{C}(Y\otimes X,Z)
\cong
\mathrm{Hom}_{C}\big(Y,Z
\otimes X^\star\big),
\end{gather*}
and similar isomorphism pulling $X$ from the second into the first component.
\end{Lemma}

\begin{proof}
Let us construct isomorphisms for the first case, all others cases are similar. We define
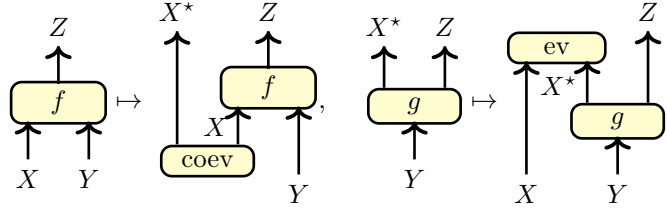
\begin{gather*}
\begin{tikzpicture}[anchorbase,scale=1]
\draw[usual,white] (0,0) to (0,0.01) node[above,black,box]{\raisebox{-0.025cm}{\hspace*{0.4cm}$f$\hspace*{0.4cm}}};
\draw[usual,directed=0.99] (-0.4,-0.45)node[below]{$X$} to (-0.4,0.05);
\draw[usual,directed=0.99] (0.4,-0.45)node[below]{$Y$} to (0.4,0.05);
\draw[usual,directed=0.99] (0,0.55) to (0,1.05)node[above]{$Z$};
\end{tikzpicture}
\mapsto
\begin{tikzpicture}[anchorbase,scale=1]
\draw[usual,white] (0,0) to (0,0.01) node[above,black,box]{\raisebox{-0.025cm}{\hspace*{0.4cm}$f$\hspace*{0.4cm}}};
\draw[usual,white] (-0.8,-0.9) to (-0.8,-0.89) node[above,black,box]{\raisebox{-0.025cm}{\hspace*{0.15cm}$\text{coev}$\hspace*{0.15cm}}};
\draw[usual,directed=0.99] (-0.4,-0.45) node[left,yshift=0.225cm]{$X$} to (-0.4,0.05);
\draw[usual,directed=0.99] (-1.2,-0.45) to (-1.2,0.05) to (-1.2,1.05)node[above]{$X^\star$};
\draw[usual,directed=0.99] (0.4,-0.85)node[below]{$Y$} to (0.4,0.05);
\draw[usual,directed=0.99] (0,0.55)to (0,1.05)node[above]{$Z$};
\end{tikzpicture}
,
\quad
\begin{tikzpicture}[anchorbase,scale=1]
\draw[usual,white] (0,0) to (0,0.01) node[above,black,box]{\raisebox{-0.025cm}{\hspace*{0.4cm}$g$\hspace*{0.4cm}}};
\draw[usual,directed=0.99] (0,-0.45)node[below]{$Y$} to (0,0.05);
\draw[usual,directed=0.99] (-0.4,0.5) to (-0.4,1.05)node[above]{$X^\star$};
\draw[usual,directed=0.99] (0.4,0.5) to (0.4,1.05)node[above]{$Z$};
\end{tikzpicture}
\mapsto
\begin{tikzpicture}[anchorbase,scale=1]
\draw[usual,white] (0,0) to (0,0.01) node[above,black,box]{\raisebox{-0.025cm}{\hspace*{0.4cm}$g$\hspace*{0.4cm}}};
\draw[usual,white] (-0.8,1) to (-0.8,1.01) node[above,black,box]{\raisebox{-0.025cm}{\hspace*{0.35cm}$\text{ev}$\hspace*{0.35cm}}};
\draw[usual,directed=0.99] (0,-0.45)node[below]{$Y$} to (0,0.05);
\draw[usual,directed=0.99] (-0.4,0.5) node[left,yshift=0.225cm]{$X^\star$} to (-0.4,1.05);
\draw[usual,directed=0.99] (-1.2,-0.45)node[below]{$X$} to (-1.2,0.55) to (-1.2,1.05);
\draw[usual,directed=0.99] (0.4,0.5) to (0.4,1.05) to (0.4,1.45)node[above]{$Z$};
\end{tikzpicture}
\end{gather*}
The fact that these are inverses follows from the zigzag relations.
\end{proof}

It thus suffices to prove that the corresponding family of algebra morphisms $f_k: \text{End}_{\mathbf{DiTL}(-2)}(\bullet^{ \otimes k}) \rightarrow \text{End}_G(V^{\otimes k})$ are all isomorphisms. By \autoref{L:tlbasis}, the dimension of $DiTL_k(-2)$ is $\frac{1}{2}\binom{2k}{k}$, which matches that of $\text{End}_G(V^{\otimes k})$ by \autoref{P:diminf}. It now suffices to show that $f_k$ is surjective for all $k$, which amounts to showing that every idempotent in $\text{End}_G(V^{\otimes k})$ is in the image of $f_k$. By \autoref{P:cgrule}, we need only show that the idempotent projecting to $V_k \subset V^{\otimes k}$ is in the image. Writing $1$ for the identity morphism on $V$, define:
\begin{gather*}
f_1 = 1, \quad  f_2 = (f_1 \otimes 1) + \frac{1}{2}(\text{ev} \circ \text{coev} + e_{\bar{1}}),
\\
f_{j+1} = (f_j \otimes 1) + (f_j \otimes 1) \circ (1^{\otimes j-1} \otimes (\text{coev} \circ \text{ev}) \otimes 1^{\otimes n-j}) \circ (f_j \otimes 1) \quad \text{for} \ j \geq 2.
\end{gather*}
Observe that this defines a sequence $\{f_j\}$ of idempotents projecting to the subrepresentation generated by the simple tensors $\{v_1^{\otimes j}, v_2^{\otimes j}\} \subset V^{\otimes j}$, which is $V_j$. Since $f_1, f_2$ are in the image of $f_1, f_2$, it follows that $f_k$ is in the image of $f_k$, and we are done.
\end{proof}

\begin{Remark}
The idempotents $\{f_j\}$ constructed above are scalar multiples of the images of the classical Jones--Wenzl projectors under restriction to the dicyclic group $\text{Dic}_{\infty} \subset SU(2).$ We elaborate on their diagrammatic interpretation in the next section. 
\end{Remark}

\begin{Theorem}\label{C:EqInf}
There is an equivalence of ribbon $\C$-linear additive categories $\tilde{\mathcal{F}}\colon\mathrm{Kar}\left( \mathbf{DiTL}(-2,\infty) \right) \to \mathbf{Rep}(\mathrm{Dic}_{\infty})$ determined by $\mathcal{F}$.
\end{Theorem}

\begin{proof}
By \autoref{T:equivalence-infinity}, the functor
\[
\mathcal{F}\colon \mathbf{DiTL}(-2,\infty)\longrightarrow \mathbf{Rep}(\mathrm{Dic}_{\infty})
\]
is a $\C$-linear additive ribbon functor which is full and faithful, and whose essential image is the full additive subcategory generated by tensor powers of $V$. By \autoref{lemma:all-appear}, every finite-dimensional simple $\mathrm{Dic}_\infty$-module occurs as a direct summand of some $V^{\otimes k}$, hence the essential image of $\mathcal{F}$ additively generates all of $\mathbf{Rep}(\mathrm{Dic}_\infty)$.

Now apply the standard universal property of the Karoubi envelope: $\mathcal{F}$ extends uniquely (up to unique isomorphism) to an additive functor
\[
\tilde{\mathcal{F}}\colon \mathrm{Kar}\!\left(\mathbf{DiTL}(-2,\infty)\right)\longrightarrow \mathbf{Rep}(\mathrm{Dic}_{\infty}),
\]
and fullness/faithfulness are preserved under idempotent completion. Since $\mathbf{Rep}(\mathrm{Dic}_\infty)$ is already idempotent complete, and every object is a direct summand of a finite direct sum of objects in the image of $\mathcal{F}$, the extended functor is essentially surjective. Therefore $\tilde{\mathcal{F}}$ is an equivalence. The ribbon structure is inherited because $\mathcal{F}$ is ribbon and the Karoubi envelope is formed inside an additive monoidal category.
\end{proof}

\subsection{Dicyclic Jones--Wenzl projectors}

We wish to generalize the above theorem to the finite degree $n$ dicyclic group $\text{Dic}_n$. In order to define the appropriate diagrammatic category for $\text{Dic}_n$, we need to start with the dicyclic Temperley--Lieb category $\mathbf{DiTL}(-2,\infty)$ and make some adjustments. First, as before, denote by $e_{\bar{1}}$ the decorated cup-cap in the dicyclic Temperley--Lieb category. 

\begin{Lemma}\label{L:idem}
Under the equivalence in \autoref{C:EqInf}, $e_{\bar{1}}$ is (a scalar multiple of) the idempotent projecting onto $\bar{\C}.$
\end{Lemma}

\begin{proof}
Under the equivalence $\tilde{\mathcal{F}}$ from \autoref{C:EqInf}, the object in $\mathrm{Kar}(\mathbf{DiTL}(-2,\infty))$ corresponding to the $2$-strand object is sent to $V^{\otimes 2}$, and $\tilde{\mathcal{F}}$ identifies endomorphism algebras:
\[
\mathrm{End}_{\mathbf{DiTL}}(2)\ \xrightarrow{\ \sim\ }\ \mathrm{End}_{\mathrm{Dic}_\infty}(V^{\otimes 2}).
\]
By definition, $e_{\bar 1}\in \mathrm{End}_{\mathrm{Dic}_\infty}(V^{\otimes 2})$ is the primitive idempotent whose image is the $\bar{\C}$-summand of $V^{\otimes 2}$. Hence its preimage under the above identification is exactly the (primitive) idempotent in $\mathrm{End}_{\mathbf{DiTL}}(2)$ cutting out the summand corresponding to $\bar{\C}$; since both sides are only defined up to the chosen normalization of the projector, this agrees up to a scalar multiple.
\end{proof}

\begin{Remark}
We will be brief below. All the classical arguments in the Temperley--Lieb case 
as e.g. in \cite[Section 3.3]{KaLi-TL-recoupling} work with the adjustment coming from \autoref{L:JWscalar} and \autoref{L:idem}. See also \cite[Section 3]{SuTuWeZh-mixed-tilting} for a detailed explanation of the corresponding ``idempotent construction strategy''.
\end{Remark}

Next, we define \emph{dicyclic Jones--Wenzl projectors}.
To this end, define $f_1 = \text{id}_{V}$ and $f_2 = f_1 + \frac{1}{2}(\text{ev} \circ \text{coev} + e_{\bar{1}})$, and let $n > 2$. Then, as in the proof of \autoref{T:equivalence-infinity}, let us recursively define a sequence of morphisms $\{f_j\}$ for $2 \leq j \leq n-1$ by:
\begin{gather}\label{Eq:JW}
f_{j+1} = f_j + f_je_{j}f_j.
\end{gather}
We slightly abuse the notation and write $f_j$ in the above expression instead of $f_j \otimes 1.$ 
We will also use the usual \emph{box notation} for these projectors.
Then, for $i\geq 2$, \autoref{Eq:JW} becomes:
\begin{gather*}
\begin{tikzpicture}[anchorbase,scale=1]
\draw[ip] (0,0) rectangle (1,0.5) node[pos=0.5]{$f_{i}$};
\end{tikzpicture}
=
\begin{tikzpicture}[anchorbase,scale=1]
\draw[ip] (0,0) rectangle (1,0.5) node[pos=0.5]{$f_{i-1}$};
\draw[usual] (1.25,0) to (1.25,0.5);
\end{tikzpicture}
+
\begin{tikzpicture}[anchorbase,scale=1]
\draw[ip] (-0.25,0) rectangle (1.5,0.5) node[pos=0.5]{$f_{i-1}$};
\draw[ip] (-0.25,0.5) rectangle (1,1) node[pos=0.5]{$f_{i-2}$};
\draw[ip] (-0.25,1) rectangle (1.5,1.5) node[pos=0.5]{$f_{i-1}$};
\draw[usual] (1.75,0) to (1.75,0.5) to[out=90,in=0] (1.5,0.7) to[out=180,in=90] (1.25,0.5);
\draw[usual] (1.75,1.5) to (1.75,1) to[out=270,in=0] (1.5,0.8) to[out=180,in=270] (1.25,1);
\end{tikzpicture}
.
\end{gather*}

\begin{Example}
The first projector $f_1$ is the identity, and for example the next two projectors in the sequence are:
\begin{gather*}
f_2 = \ \begin{tikzpicture}[anchorbase, scale=0.5]
\draw [draw=none] (1,0) rectangle (2,2);
\draw [usual] (1,0) [dnup=0];
\draw [usual] (2,0) [dnup=0];
\end{tikzpicture}
+ \frac{1}{2} \ 
\begin{tikzpicture}[anchorbase, scale=0.5]
\draw [draw=none] (0,0) rectangle (1,2);
\draw [usual] (0,0) [dndn=1];
\draw [usual] (0,2) [upup=1];
\end{tikzpicture}
+ \frac{1}{2} \ \begin{tikzpicture}[anchorbase, scale=0.5]
\draw [draw=none] (0,0) rectangle (2,2);
\draw [usual] [postaction={decorate,decoration={markings,mark=at position 0.5 with {\node[circle,fill=blue,inner sep=1.5pt] {};}}}] (0,0) [dndn=1];
\draw [usual] [postaction={decorate,decoration={markings,mark=at position 0.5 with {\node[circle,fill=blue,inner sep=1.5pt] {};}}}] (0,2) [upup=1];
\end{tikzpicture},
\\
\begin{aligned}
f_3 =& \ \begin{tikzpicture}[anchorbase, scale=0.5]
\draw [draw=none] (1,0) rectangle (2,2);
\draw [usual] (1,0) [dnup=0];
\draw [usual](2,0) [dnup=0];
\draw [usual](3,0) [dnup=0];
\end{tikzpicture}
\ + \
\begin{tikzpicture}[anchorbase, scale=0.5]
\draw [draw=none] (1,0) rectangle (2,2);
\draw [usual] (1,0) [dnup=0];
\draw [usual](2,0) [dndn=1];
\draw [usual](2,2) [upup=1];
\end{tikzpicture}
+  \ \frac{3}{4} \left( \  \begin{tikzpicture}[anchorbase, scale=0.5]
\draw [draw=none] (1,0) rectangle (2,2);
\draw [usual](1,0) [dndn=1];
\draw [usual](1,2) [upup=1];
\draw [usual](3,0) [dnup=0];
\end{tikzpicture}
\ + \ \begin{tikzpicture}[anchorbase, scale=0.5]
\draw [draw=none] (1,0) rectangle (2,2);
\draw [usual][postaction={decorate,decoration={markings,mark=at position 0.5 with {\node[circle,fill=blue,inner sep=1.5pt] {};}}}] (1,0) [dndn=1];
\draw [usual][postaction={decorate,decoration={markings,mark=at position 0.5 with {\node[circle,fill=blue,inner sep=1.5pt] {};}}}] (1,2) [upup=1];
\draw [usual](3,0) [dnup=0];
\end{tikzpicture} \ \right)
\ + \ \frac{1}{4} \left( \ 
\begin{tikzpicture}[anchorbase, scale=0.5]
\draw [draw=none] (1,0) rectangle (2,2);
\draw [usual](1,0) [dndn=1];
\draw [usual][postaction={decorate,decoration={markings,mark=at position 0.5 with {\node[circle,fill=blue,inner sep=1.5pt] {};}}}] (1,2) [upup=1];
\draw [usual][postaction={decorate,decoration={markings,mark=at position 0.5 with {\node[circle,fill=blue,inner sep=1.5pt] {};}}}] (3,0) [dnup=0];
\end{tikzpicture}
\ + \ 
\begin{tikzpicture}[anchorbase, scale=0.5]
\draw [draw=none] (1,0) rectangle (2,2);
\draw [usual][postaction={decorate,decoration={markings,mark=at position 0.5 with {\node[circle,fill=blue,inner sep=1.5pt] {};}}}] (1,0) [dndn=1];
\draw [usual](1,2) [upup=1];
\draw [usual][postaction={decorate,decoration={markings,mark=at position 0.5 with {\node[circle,fill=blue,inner sep=1.5pt] {};}}}] (3,0) [dnup=0];
\end{tikzpicture} \ \right) 
\\
&+ \ \frac{1}{2} \left( \ \begin{tikzpicture}[anchorbase, scale=0.5]
\draw [draw=none] (1,0) rectangle (2,2);
\draw [usual](1,0) [dndn=1];
\draw [usual](3,0) [dnup=-2];
\draw [usual](2,2) [upup=1];
\end{tikzpicture}
\ + \ \begin{tikzpicture}[anchorbase, scale=0.5]
\draw [draw=none] (1,0) rectangle (2,2);
\draw [usual][postaction={decorate,decoration={markings,mark=at position 0.5 with {\node[circle,fill=blue,inner sep=1.5pt] {};}}}] (1,0) [dndn=1];
\draw [usual][postaction={decorate,decoration={markings,mark=at position 0.5 with {\node[circle,fill=blue,inner sep=1.5pt] {};}}}] (3,0) [dnup=-2];
\draw [usual](2,2) [upup=1];
\end{tikzpicture}
\ + \ \begin{tikzpicture}[anchorbase, scale=0.5]
\draw [draw=none] (1,0) rectangle (2,2);
\draw [usual](1,2) [upup=1];
\draw [usual](1,0) [dnup=2];
\draw [usual](2,0) [dndn=1];
\end{tikzpicture}
\ + \ \begin{tikzpicture}[anchorbase, scale=0.5]
\draw [draw=none] (1,0) rectangle (2,2);
\draw [usual][postaction={decorate,decoration={markings,mark=at position 0.5 with {\node[circle,fill=blue,inner sep=1.5pt] {};}}}] (1,2) [upup=1];
\draw [usual][postaction={decorate,decoration={markings,mark=at position 0.5 with {\node[circle,fill=blue,inner sep=1.5pt] {};}}}] (1,0) [dnup=2];
\draw [usual](2,0) [dndn=1];
\end{tikzpicture} \ \right).
\end{aligned}
\end{gather*}
This follows by direct calculation.
\end{Example}

We state the following properties (written in terms of the category but omitting $\otimes\text{id}$ etc.), where the indexes are all admissible possibilities:
\begin{gather}\label{Eq:JWequations}
f_{i}\neq 0,
\quad
f_{i}f_{i}=f_{i},
\quad
f_{i}\text{coev}=\text{ev}f_{i}=0,
\quad
f_{i}e_{\bar{1}}=e_{\bar{1}}f_{i}=0.
\end{gather}
In box notation the relations take a familiar form:
\begin{gather*}
\begin{tikzpicture}[anchorbase,scale=1]
\draw[ip] (0,0) rectangle (1,0.5) node[pos=0.5]{$f_{i}$};
\draw[ip] (0,0.5) rectangle (1,1) node[pos=0.5]{$f_{i}$};
\end{tikzpicture}
=
\begin{tikzpicture}[anchorbase,scale=1]
\draw[ip] (0,0) rectangle (1,0.5) node[pos=0.5]{$f_{i}$};
\end{tikzpicture}
,\quad
\begin{tikzpicture}[anchorbase,scale=1]
\draw[usual,white] (0.25,0.5) to[out=90,in=180] (0.5,0.75) to[out=0,in=90] (0.75,0.5);
\draw[ip] (0,0) rectangle (1,0.5) node[pos=0.5]{$f_{i}$};
\draw[usual] (0.25,0) to[out=270,in=180] (0.5,-0.25) to[out=0,in=270] (0.75,0);
\end{tikzpicture}
=
0
=
\begin{tikzpicture}[anchorbase,scale=1]
\draw[usual,white] (0.25,0) to[out=270,in=180] (0.5,-0.25) to[out=0,in=270] (0.75,0);
\draw[ip] (0,0) rectangle (1,0.5) node[pos=0.5]{$f_{i}$};
\draw[usual] (0.25,0.5) to[out=90,in=180] (0.5,0.75) to[out=0,in=90] (0.75,0.5);
\end{tikzpicture}
,\quad
\begin{tikzpicture}[anchorbase,scale=1]
\draw[usual,white] (0.25,0.5) to[out=90,in=180] (0.5,0.75) to[out=0,in=90] (0.75,0.5);
\draw[ip] (0,0) rectangle (1,0.5) node[pos=0.5]{$f_{i}$};
\draw[usual][postaction={decorate,decoration={markings,mark=at position 0.5 with {\node[circle,fill=blue,inner sep=1.5pt] {};}}}] (0.25,0) to[out=270,in=180] (0.5,-0.25) to[out=0,in=270] (0.75,0);
\end{tikzpicture}
=
0
=
\begin{tikzpicture}[anchorbase,scale=1]
\draw[usual,white] (0.25,0) to[out=270,in=180] (0.5,-0.25) to[out=0,in=270] (0.75,0);
\draw[ip] (0,0) rectangle (1,0.5) node[pos=0.5]{$f_{i}$};
\draw[usual][postaction={decorate,decoration={markings,mark=at position 0.5 with {\node[circle,fill=blue,inner sep=1.5pt] {};}}}] (0.25,0.5) to[out=90,in=180] (0.5,0.75) to[out=0,in=90] (0.75,0.5);
\end{tikzpicture}
.
\end{gather*}

\begin{Lemma}
\autoref{Eq:JWequations} hold and characterize the dicyclic Jones--Wenzl projectors uniquely.
\end{Lemma}

\begin{proof}
Existence is by construction: by \autoref{L:JWscalar} the chosen normalization differs from the usual Jones--Wenzl recursion only by a fixed nonzero scalar, so the standard inductive construction produces idempotents satisfying \autoref{Eq:JWequations}.

For uniqueness, let $p_n$ and $p_n'$ be two solutions of \autoref{Eq:JWequations} in the relevant endomorphism algebra. Then $d=p_n-p_n'$ satisfies the same annihilation relations as $p_n$ (and $p_n'$), so it factors through the sum of the ``lower'' through-degree summands. Equivalently, $d$ lies in the ideal generated by the $e_i$ (and, in the dicyclic case, $e_{\bar 1}$) that are killed by \autoref{Eq:JWequations}. But \autoref{Eq:JWequations} also force $p_n d=d p_n'=d$ and $d^2=0$, hence $d=0$. Thus $p_n=p_n'$, so the equations characterize the projector uniquely.
\end{proof}

\begin{Lemma}
The dicyclic Jones--Wenzl projectors satisfy the following.
\begin{enumerate}

\item We have \emph{hom vanishing}, i.e. for $0\leq j\leq i$ we have
\begin{gather}\label{10-eq:jwhom}
f_{j}
\mathrm{Hom}_{\mathbf{DiTL}(-2,\infty)}
(\bullet^{i},\bullet^{j})f_{i}
=
\begin{cases*}
\C\{f_{i}\}&\text{if }i=j,
\\
0&\text{else}.
\end{cases*}
\end{gather}
Similarly if $j\geq i$.

\item We have \emph{absorption}, i.e.
\begin{gather*}
\begin{tikzpicture}[anchorbase,scale=1]
\draw[ip] (-0.25,0) rectangle (1.25,0.5) node[pos=0.5]{$f_{i}$};
\draw[ip] (0,0.5) rectangle (1,1) node[pos=0.5]{$f_{j}$};
\draw[usual] (-0.15,0.5) to (-0.15,1);
\draw[usual] (1.15,0.5) to (1.15,1);
\end{tikzpicture}
=
\begin{tikzpicture}[anchorbase,scale=1]
\draw[ip] (0,0) rectangle (1,0.5) node[pos=0.5]{$f_{i}$};
\end{tikzpicture}
=
\begin{tikzpicture}[anchorbase,scale=1]
\draw[ip] (-0.25,0) rectangle (1.25,0.5) node[pos=0.5]{$f_{i}$};
\draw[ip] (0,-0.5) rectangle (1,0) node[pos=0.5]{$f_{j}$};
\draw[usual] (-0.15,-0.5) to (-0.15,0);
\draw[usual] (1.15,-0.5) to (1.15,0);
\end{tikzpicture}
\quad\text{where }0\leq j\leq i.
\end{gather*}

\item We have \emph{partial trace properties}, i.e. for $i\geq 2$:
\begin{gather*}
\begin{tikzpicture}[anchorbase,scale=1]
\draw[ip] (-0.5,0) rectangle (0.5,0.5) node[pos=0.5]{$f_{i}$};
\draw[usual] (0.25,0.5) to[out=90,in=180] (0.5,0.75) to[out=0,in=90] (0.75,0.5) to (0.75,0) to[out=270,in=0] (0.5,-0.25) to[out=180,in=270] (0.25,0);
\end{tikzpicture}
=
\begin{tikzpicture}[anchorbase,scale=1]
\draw[ip] (-0.5,0) rectangle (0.5,0.5) node[pos=0.5]{$f_{i}$};
\draw[usual][postaction={decorate,decoration={markings,mark=at position 0.2 with {\node[circle,fill=blue,inner sep=1.5pt] {};}}}][postaction={decorate,decoration={markings,mark=at position 0.8 with {\node[circle,fill=blue,inner sep=1.5pt] {};}}}] (0.25,0.5) to[out=90,in=180] (0.5,0.75) to[out=0,in=90] (0.75,0.5) to (0.75,0) to[out=270,in=0] (0.5,-0.25) to[out=180,in=270] (0.25,0);
\end{tikzpicture}
=
-
\begin{tikzpicture}[anchorbase,scale=1]
\draw[ip] (0,0) rectangle (1,0.5) node[pos=0.5]{$f_{i-1}$};
\end{tikzpicture}
.
\end{gather*}
Similar partial traces, but with different scalars, hold for $i=0,1$.

\end{enumerate}
\end{Lemma}

\begin{proof}
One can copy the calculations in \cite[Section 10]{Tu-qt}.
\end{proof}

The following justifies the name dicyclic Jones--Wenzl projectors.

\begin{Proposition}
Under the equivalence with $\mathbf{Rep}(\text{Dic}_{\infty})$, the 
dicyclic Jones--Wenzl projectors 
are the idempotents projecting to the leading summand of $V^{\otimes k}$.
This also holds for the finite group for $1 \leq k <n$. 
\end{Proposition}

\begin{proof}
Let $p_k$ be the $k$-strand dicyclic Jones--Wenzl projector in
$\mathrm{End}_{\mathbf{DiTL}(-2,\infty)}(k)$, and transport it along the equivalence
$\tilde{\mathcal F}$ of \autoref{C:EqInf} to an idempotent
$\tilde{\mathcal F}(p_k)\in \mathrm{End}_G(V^{\otimes k})$.

By \autoref{Eq:JWequations}, the projector $p_k$ is characterized by being an
idempotent which is annihilated by the ``lowering'' generators (the $e_i$,
and in the dicyclic case also the $e_{\bar 1}$ at the appropriate boundary),
and by having the usual normalization (up to the scalar of \autoref{L:JWscalar}).
Under $\tilde{\mathcal F}$, these generators correspond to the canonical
$G$-intertwiners cutting out the non-leading summands (cf. the decomposition of
$V^{\otimes 2}$ and \autoref{L:idem}), so the condition that they annihilate
$\tilde{\mathcal F}(p_k)$ is exactly the condition that $\tilde{\mathcal F}(p_k)$
kills every summand of $V^{\otimes k}$ that factors through a smaller tensor power.
Hence $\tilde{\mathcal F}(p_k)$ acts as the identity on the unique summand
of $V^{\otimes k}$ of maximal through-degree (the ``leading'' summand) and as
zero on all other summands.
By the uniqueness part of the Jones--Wenzl characterization,
this determines $\tilde{\mathcal F}(p_k)$ as the idempotent projector onto that
leading summand.

For $\mathrm{Dic}_n$, the same argument applies for $1\le k<n$: in this range the
diagrammatic relations (and hence the Jones--Wenzl recursion/characterization)
agree with the infinite case, so the transported projector is still the
projection onto the leading summand of $V^{\otimes k}$.
\end{proof}

\begin{Remark}
The recursive formula \autoref{Eq:JW} for these idempotents is the dicyclic version of the Jones--Wenzl idempotent recursion with one crucial difference: the appearing scalar is $1$ in \autoref{Eq:JW} and $k/(k+1)$ for the Temperley--Lieb algebra. The difference comes from 
\autoref{L:JWscalar} where the ratio would be $k/(k+1)$ for $SU(2)$.
\end{Remark}

\section{Finite dicyclic categories}\label{S:Finite}

We now study the finite case.

\subsection{The degree $n$ dicyclic Temperley--Lieb category $\mathbf{DiTL}(\delta,n)$}

Let us write $e_i$ short for either the evaluation or the coevaluation applied in spot $i$.
We have the following analog of $\mathbf{DiTL}(\delta, \infty)$:

\begin{Definition}
The \emph{dicyclic Temperley--Lieb category of degree $n$} $\mathbf{DiTL}(\delta, n)$
is defined exactly as $\mathbf{DiTL}(\delta, \infty)$ but with the additional morphism
$J=J_n \in \text{End}(\bullet^{\otimes n})$, the \emph{degree $n$ dicyclic Jones--Wenzl projector}, and the additional relations:
\begin{gather*}
(1) \ f_jJ = J = Jf_j, \ 1 \leq j \leq n-1,\quad \quad
(2) \ J^2 = J, \quad \quad (3) \ f_{n}J = 0 = Jf_{n}, \\
(4) \ e_nJe_n = -\frac{1}{2}f_{n-1},
\quad \quad
(5) \ e_iJ = 0 = Je_i, \ 1 \leq i < n, \quad \quad (6) \ e_iJ = Je_i, \ i > n, \quad \quad (7) \ J\otimes\text{id} = -2Je_nJ,
\end{gather*}
where the dicyclic Jones--Wenzl projectors $\{f_j\} \in \text{End}(\bullet^{\otimes j})$ for $1 \leq j \leq n$ are defined as before.
\end{Definition}

As alluded to above, the relations in the definition of $\mathbf{DiTL}(\delta,n)$ are type D analogs of the relations satisfied by the classical Jones--Wenzl projectors.

\begin{Lemma}\label{L:ribbon2}
The category $\mathbf{DiTL}(-2,n)$ is ribbon by inheriting (as an extension) the structures from 
$\mathbf{DiTL}(\delta, \infty)$.
\end{Lemma}

\begin{proof}
As in \autoref{L:ribbon}.
\end{proof}

We, inspired by \cite{Co-jelly}, would like to draw a jellyfish for $J$. Or, alternatively, an $2n$ valent vertex, a notation we learned from \cite{El-two-color-soergel}. 
However, these are taken, so we instead draw:
\begin{gather*}
J\leftrightsquigarrow
\begin{tikzpicture}[anchorbase,scale=1]
\draw[ipp] (0,0) rectangle (1,0.5) node[pos=0.5]{$J$};
\end{tikzpicture}
.
\end{gather*}
The relations then take the following form:
\begin{gather*}
\begin{tikzpicture}[anchorbase,scale=1]
\draw[ipp] (-0.25,0) rectangle (1.25,0.5) node[pos=0.5]{$J$};
\draw[ip] (0,0.5) rectangle (1,1) node[pos=0.5]{$f_{j}$};
\draw[usual] (-0.15,0.5) to (-0.15,1);
\draw[usual] (1.15,0.5) to (1.15,1);
\end{tikzpicture}
=
\begin{tikzpicture}[anchorbase,scale=1]
\draw[ipp] (0,0) rectangle (1,0.5) node[pos=0.5]{$J$};
\end{tikzpicture}
=
\begin{tikzpicture}[anchorbase,scale=1]
\draw[ipp] (-0.25,0) rectangle (1.25,0.5) node[pos=0.5]{$J$};
\draw[ip] (0,-0.5) rectangle (1,0) node[pos=0.5]{$f_{j}$};
\draw[usual] (-0.15,-0.5) to (-0.15,0);
\draw[usual] (1.15,-0.5) to (1.15,0);
\end{tikzpicture}
\quad\text{where }0\leq j\leq i,
\quad
\begin{tikzpicture}[anchorbase,scale=1]
\draw[ipp] (0,0) rectangle (1,0.5) node[pos=0.5]{$J$};
\draw[ipp] (0,0.5) rectangle (1,1) node[pos=0.5]{$J$};
\end{tikzpicture}
=
\begin{tikzpicture}[anchorbase,scale=1]
\draw[ipp] (0,0) rectangle (1,0.5) node[pos=0.5]{$J$};
\end{tikzpicture}
,\quad
\begin{tikzpicture}[anchorbase,scale=1]
\draw[ip] (0,0) rectangle (1,0.5) node[pos=0.5]{$f_{n}$};
\draw[ipp] (0,0.5) rectangle (1,1) node[pos=0.5]{$J$};
\end{tikzpicture}
=
0
=
\begin{tikzpicture}[anchorbase,scale=1]
\draw[ipp] (0,0) rectangle (1,0.5) node[pos=0.5]{$J$};
\draw[ip] (0,0.5) rectangle (1,1) node[pos=0.5]{$f_{n}$};
\end{tikzpicture}
,
\end{gather*}
for (1) to (3). The relations (6) is a far commutativity relation and not displayed, while
\begin{gather*}
\begin{tikzpicture}[anchorbase,scale=1]
\draw[ipp] (-0.5,0) rectangle (0.5,0.5) node[pos=0.5]{$J$};
\draw[usual] (0.25,0.5) to[out=90,in=180] (0.5,0.75) to[out=0,in=90] (0.75,0.5) to (0.75,0) to[out=270,in=0] (0.5,-0.25) to[out=180,in=270] (0.25,0);
\end{tikzpicture}
=
-\frac{1}{2}\cdot
\begin{tikzpicture}[anchorbase,scale=1]
\draw[ip] (0,0) rectangle (1,0.5) node[pos=0.5]{$f_{n-1}$};
\end{tikzpicture}
,\quad
\begin{tikzpicture}[anchorbase,scale=1]
\draw[usual,white] (0.25,0.5) to[out=90,in=180] (0.5,0.75) to[out=0,in=90] (0.75,0.5);
\draw[ipp] (0,0) rectangle (1,0.5) node[pos=0.5]{$J$};
\draw[usual] (0.25,0) to[out=270,in=180] (0.5,-0.25) to[out=0,in=270] (0.75,0);
\end{tikzpicture}
=
0
=
\begin{tikzpicture}[anchorbase,scale=1]
\draw[usual,white] (0.25,0) to[out=270,in=180] (0.5,-0.25) to[out=0,in=270] (0.75,0);
\draw[ipp] (0,0) rectangle (1,0.5) node[pos=0.5]{$J$};
\draw[usual] (0.25,0.5) to[out=90,in=180] (0.5,0.75) to[out=0,in=90] (0.75,0.5);
\end{tikzpicture}
,
\end{gather*}
are relations (4) and (5). The final relation (7) is (equivalent to):
\begin{gather*}
\begin{tikzpicture}[anchorbase,scale=1]
\draw[ipp] (0,0) rectangle (1,0.5) node[pos=0.5]{$J$};
\draw[usual] (1.25,0) to (1.25,0.5);
\end{tikzpicture}
=-2\cdot
\begin{tikzpicture}[anchorbase,scale=1]
\draw[ipp] (-0.25,0) rectangle (1.5,0.5) node[pos=0.5]{$J$};
\draw[ip] (-0.25,0.5) rectangle (1,1) node[pos=0.5]{$f_{i-2}$};
\draw[ipp] (-0.25,1) rectangle (1.5,1.5) node[pos=0.5]{$J$};
\draw[usual] (1.75,0) to (1.75,0.5) to[out=90,in=0] (1.5,0.7) to[out=180,in=90] (1.25,0.5);
\draw[usual] (1.75,1.5) to (1.75,1) to[out=270,in=0] (1.5,0.8) to[out=180,in=270] (1.25,1);
\end{tikzpicture}
.
\end{gather*}

\begin{Remark}
Note that
\begin{gather*}
\begin{tikzpicture}[anchorbase,scale=1]
\draw[ipp] (0,0) rectangle (1,0.5) node[pos=0.5]{$J$};
\end{tikzpicture}
\neq
\begin{tikzpicture}[anchorbase,scale=1]
\draw[ip] (0,0) rectangle (1,0.5) node[pos=0.5]{$f_{n}$};
\end{tikzpicture}
,
\end{gather*}
and indeed their partial trace scalar down to $f_{n-1}$ is different. Moreover, relation (7) prevents a new idempotent of the same type as $J$ to appear.
\end{Remark}

\subsection{An equivalence $\text{Kar}({\mathbf{DiTL}}(-2, n) \cong \mathbf{Rep}(\text{Dic}_n)$}

We are ready to prove our final statements.
Let $\mathbf{Rep}(\text{Dic}_{n})$ denote the category of finite dimensional complex representations of $\text{Dic}_{n}$, and let $\mathcal{C}_n$ denote the full monoidal subcategory of $\mathbf{Rep}(\text{Dic}_{n})$ monoidally generated by $V$.

\begin{Lemma}\label{L:functorTwo}
We have the following.
\begin{enumerate}

\item The morphisms $\text{coev}$, $\text{ev}$, and $e_{\bar{1}}$ monoidally generate $\mathcal{C}_n$ together with the projector $J_{rt}$ factoring $V^{\otimes n}\to V_{n'}\to V^{\otimes n}$.

\item The three categories $\mathbf{DiTL}(-2,n)$, $\mathcal{C}_n$ and $\mathbf{Rep}(\text{Dic}_{n})$, are ribbon with structures coming from their infinite counterparts.

\item There is a well-defined ribbon $\C$-linear functor $\mathcal{F}_n\colon\mathbf{DiTL}(-2,n) \to\mathcal{C}_n$ determined by mapping $\bullet \mapsto V$ and the cap, cup, and decorated cup-cap morphisms to $\text{ev}, \text{coev}, $ and $e_{\bar{1}}$ respectively, and $J$ to $J_{rt}$.

\end{enumerate}
\end{Lemma}

\begin{proof}
Using \autoref{P:fusion-graph}, this is
not much different from the infinite case. The only interesting thing to check are the relations involving $J$. However, all of them follow directly from \autoref{P:fusion-graph} and the usual properties of Jones--Wenzl type projectors.
\end{proof}

\begin{Theorem}\label{T:equivalence-infinityTwo}
The functor $\mathcal{F}_n$ from \autoref{L:functorTwo} is fully faithful.
\end{Theorem}

\begin{proof}
Note that the generating morphism $J$ can only exist in hom sets $\text{Hom}(\bullet^{\otimes k}, \bullet^ {\otimes m})$ where one of $k,m \geq n.$ We have established that the dicyclic Temperley--Lieb category $\mathbf{DiTL}(-2, \infty)$ corresponds to the representation theory of the infinite dicyclic group $\text{Dic}_{\infty}$. The additional generating morphism $J$ in the degree $n$ category $\mathbf{DiTL}(\delta,n)$ corresponds to the idempotent obtained by passing to the finite group $\text{Dic}_n$ and studying the breaking point in the decomposition of $V^{\otimes n}$; this is the idempotent projecting to $V_{n'} \subset V^{\otimes n}$ discussed previously in \autoref{P:fusion-graph}. 

The defining relations for $f_1,f_2,\dots,f_n,J$ are satisfied by the Jones--Wenzl idempotents $\text{Res}^{SU(2)}_{Dic_n}(JW_k)$, see \autoref{L:functorTwo}. To see that $\tilde{\mathcal{F}}$ is fully faithful one can, as in the infinite case, restrict to endomorphism algebras. Then \cite[Proposition 1.50]{BaBeHa-mckay} applies.
\end{proof}

\begin{Theorem}\label{C:EqInfTwo}
There is an equivalence of ribbon $\C$-linear additive categories $\tilde{\mathcal{F}}_n\colon\mathrm{Kar}\left( \mathbf{DiTL}(-2,n) \right) \to\mathbf{Rep}(\mathrm{Dic}_{n})$ determined by $\mathcal{F}_n$.
\end{Theorem}

\begin{proof}
By \autoref{T:equivalence-infinityTwo}, the functor
\[
\mathcal{F}_n\colon \mathbf{DiTL}(-2,n)\longrightarrow \mathbf{Rep}(\mathrm{Dic}_{n})
\]
is a $\C$-linear additive ribbon functor which is full and faithful, and whose
essential image is the full additive subcategory generated by the tensor powers
$V^{\otimes k}$ (with $1\le k\le n-1$, and their direct sums/summands).

Passing to the Karoubi envelope, the universal property gives an induced
$\C$-linear additive ribbon functor
\[
\tilde{\mathcal{F}}_n\colon \mathrm{Kar}\!\left(\mathbf{DiTL}(-2,n)\right)\longrightarrow
\mathbf{Rep}(\mathrm{Dic}_{n}),
\]
and fullness/faithfulness are preserved under idempotent completion. Since
$\mathbf{Rep}(\mathrm{Dic}_n)$ is idempotent complete and, by
\autoref{lemma:all-appear} (applied to the finite quotient), every simple
$\mathrm{Dic}_n$-module occurs as a direct summand of some $V^{\otimes k}$,
the image of $\tilde{\mathcal{F}}_n$ is essentially surjective. Hence,
$\tilde{\mathcal{F}}_n$ is an equivalence of ribbon $\C$-linear additive
categories.
\end{proof}

\section{Dicyclic bases and walks}

\addtocounter{subsection}{1}

We explain how to construct a basis element for the centralizer algebra $\text{End}_G(V^{\otimes k})$ where $G = \text{Dic}_{n}$, or equivalently the dicyclic Temperley--Lieb algebra $DiTL_k(-2,n)$, given a closed walk on the type $D$ Dynkin diagram. We first briefly recall this well known construction in the classical case of $G = SU(2)$ (equivalently $SL(2,\mathbb{C})$ or $\mathfrak{sl}_2(\mathbb{C})$), which can be found in many places (see e.g. \cite{Tu-qt}). A \emph{closed walk} of length $2k$ on the type $A_{\infty}$ diagram
\begin{gather*}
A_{\infty}=
\dynkin[scale=2]{A}{oooo}
\;\cdots
\leftrightsquigarrow
\N
\end{gather*}
is a sequence $w \subset \{+1, -1\}^{2k}$ (with symbols written $+,-$) such that $\sum_{i=1}^m w_i \geq 0$ for all $1 \leq m < 2k$ and $\sum_{i=1}^{2k} w_i = 0.$ For example, $w = (+, +, -, +, -, -)$ is a closed walk of length 6. The number of such walks is well known to be the Catalan number $C_n$. The classical construction of Temperley\textendash Lieb diagrams from a closed type $A$ walk $w$ of length $2k$ is as follows. For every $+ \in w$, one draws a strand, and for every $- \in w$, the nearest strand is ``cupped off'' and the final diagram is bent over. 
For example, let us take the length 6 walk $w = (+, +, -, +, -, -) $ as before. We will produce a 3-strand TL diagram. We first draw two strands, capping off the second (for the first three steps $+$, $+$, $-$). Then following this procedure for the rest of the sequence yields:
\begin{gather*}
\emptyset\quad \rightarrow \quad
\begin{tikzpicture}[anchorbase, scale=0.5]
\draw [draw=none] (1,0) rectangle (2,2);
\draw [usual] (1,0) [dnup=0];
\end{tikzpicture} \quad \rightarrow \quad
\begin{tikzpicture}[anchorbase, scale=0.5]
\draw [draw=none] (1,0) rectangle (2,2);
\draw [usual] (1,0) [dnup=0];
\draw [usual] (2,0) [dnup=0];
\end{tikzpicture} \quad \rightarrow \quad \begin{tikzpicture}[anchorbase, scale=0.5]
\draw [draw=none] (1,0) rectangle (2,2);
\draw [usual] (1,0) [dnup=0];
\draw [usual] (2,2) [upup=1];
\end{tikzpicture} \quad \rightarrow \quad 
\begin{tikzpicture}[anchorbase, scale=0.5]
\draw [draw=none] (1,0) rectangle (4,2);
\draw [usual] (1,0) [dnup=0];
\draw [usual] (2,2) [upup=1];
\draw [usual] (4,0) [dnup=0];
\end{tikzpicture} \quad \rightarrow \quad
\begin{tikzpicture}[anchorbase, scale=0.5]
\draw [draw=none] (1,0) rectangle (4,2);
\draw [usual] (1,0) [dnup=0];
\draw [usual] (2,2) [upup=1];
\draw [usual] (4,2) [upup=1];
\end{tikzpicture} \quad \rightarrow \quad
\begin{tikzpicture}[anchorbase, scale=0.5]
\draw [draw=none] (1,0) rectangle (6,2);
\draw[usual]
  (1,2) to[out=-90,in=-90,looseness=1.0] (6,2);
\draw [usual] (2,2) [upup=1];
\draw [usual] (4,2) [upup=1];
\end{tikzpicture}
.
\end{gather*}
This procedure gives a diagram $\varnothing \to 2k$, which we bend over, bending the $k$ rightmost strands, to get a diagram $k\to k$. For example:
\begin{gather*}
\begin{tikzpicture}[anchorbase, scale=0.5]
\draw [draw=none] (1,0) rectangle (6,2);
\draw[usual]
(1,2) to[out=-90,in=-90,looseness=1.0] (6,2);
\draw [usual] (2,2) [upup=1];
\draw [usual] (4,2) [upup=1];
\draw [orchid] (6,2) [dndn=1] to (7,0);
\draw[orchid]
(5,2) to[out=90,in=90,looseness=1.1] (8,2) to (8,0);
\draw[orchid]
(4,2) to[out=90,in=90,looseness=1.0] (9,2) to (9,0);
\end{tikzpicture}
\xrightarrow{\text{bend}}
\begin{tikzpicture}[anchorbase, scale=0.5]
\draw [draw=none] (1,0) rectangle (2,2);
\draw [usual] (1,0) [dnup=0];
\draw [usual](2,0) [dndn=1];
\draw [usual](2,2) [upup=1];
\end{tikzpicture}
.
\end{gather*}
The resulting diagrams form the usual noncrossing partition basis of the Temperley--Lieb algebra. By bending \autoref{L:Bending}, this procedure gives bases for all hom-spaces in the Temperley--Lieb category as well.

Using walks on
\begin{gather*}
\tilde{D}_{n}=
\dynkin[
labels={0',0,1,2,,n{-}2,n{-}1,n',n},
label directions={,,left,,,,right,,},
scale=1.8,
extended
] D{oooo...oo*o}.
\end{gather*}
we now describe a basis of $DiTL_k(-2,n)$. For $DiTL_k(-2)$ the same construction works by taking $n\gg k$, so we henceforth only describe the situation for $DiTL_k(-2,n)$.

A \emph{closed walk} of length $2k$ on $DiTL_k(-2,n)$ is a walk starting at $0$ and ending at $0$. We encode such a walk as a sequence $w\subset\{+,-,\bplus,\bminus,\rplus,\rminus\}^{2k}$ where $\pm$ step from $i$ to $i\pm 1$, and $\bplus$ from $0'$ to $1$, $\bminus$ from $1$ to $0'$, and $\rplus$ from $n-1$ to $n'$, $\rminus$ from $n'$ to $n-1$.

The diagrammatic procedure is now as follows. Every $\bplus$ draws a strand with a blue dot, and $\bminus$ cups off a strand adding a blue dot. For example, $(+,\bminus,\bplus,-)$ is the walk $0\to 1\to 0'\to 1\to 0$ and gives
\begin{gather*}
\emptyset\quad \rightarrow \quad
\begin{tikzpicture}[anchorbase, scale=0.5]
\draw [draw=none] (1,0) rectangle (2,2);
\draw [usual] (1,0) [dnup=0];
\end{tikzpicture} 
\quad \rightarrow \quad
\begin{tikzpicture}[anchorbase, scale=0.5]
\draw [draw=none] (1,0) rectangle (2,2);
\draw[usual][postaction={decorate,decoration={markings,mark=at position 0.5 with {\node[circle,fill=blue,inner sep=1.5pt] {};}}}] (1,2) [upup=1];
\end{tikzpicture} 
\quad \rightarrow \quad 
\begin{tikzpicture}[anchorbase, scale=0.5]
\draw [draw=none] (1,0) rectangle (2,2);
\draw[usual][postaction={decorate,decoration={markings,mark=at position 0.5 with {\node[circle,fill=blue,inner sep=1.5pt] {};}}}] (1,2) [upup=1];
\draw [usual][postaction={decorate,decoration={markings,mark=at position 0.5 with {\node[circle,fill=blue,inner sep=1.5pt] {};}}}] (3,0) [dnup=0];
\end{tikzpicture} 
\quad \rightarrow \quad 
\begin{tikzpicture}[anchorbase, scale=0.5]
\draw [draw=none] (1,0) rectangle (2,2);
\draw[usual][postaction={decorate,decoration={markings,mark=at position 0.5 with {\node[circle,fill=blue,inner sep=1.5pt] {};}}}] (1,2) [upup=1];
\draw[usual][postaction={decorate,decoration={markings,mark=at position 0.5 with {\node[circle,fill=blue,inner sep=1.5pt] {};}}}] (3,2) [upup=1];
\end{tikzpicture} 
\xrightarrow{\text{bend}}
\begin{tikzpicture}[anchorbase, scale=0.5]
\draw [draw=none] (1,0) rectangle (2,2);
\draw [usual][postaction={decorate,decoration={markings,mark=at position 0.5 with {\node[circle,fill=blue,inner sep=1.5pt] {};}}}](1,0) [dndn=1];
\draw [usual][postaction={decorate,decoration={markings,mark=at position 0.5 with {\node[circle,fill=blue,inner sep=1.5pt] {};}}}](1,2) [upup=1];
\end{tikzpicture}
.
\end{gather*}
For $\rplus$ and $\rminus$ we do the same but with red dots and reversed roles of plus and minus, e.g. $(\dots,+,\rplus,\rminus,-,\dots)$ as part of a walk corresponds to
\begin{gather*}
\begin{tikzpicture}[anchorbase, scale=0.5]
\draw [draw=none] (1,0) rectangle (2,2);
\draw [usual] (1,0) [dnup=0];
\end{tikzpicture} 
\quad \rightarrow \quad
\begin{tikzpicture}[anchorbase, scale=0.5]
\draw [draw=none] (1,0) rectangle (2,2);
\draw[usual][postaction={decorate,decoration={markings,mark=at position 0.5 with {\node[circle,fill=red,inner sep=1.5pt] {};}}}] (1,2) [upup=1];
\end{tikzpicture} 
\quad \rightarrow \quad 
\begin{tikzpicture}[anchorbase, scale=0.5]
\draw [draw=none] (1,0) rectangle (2,2);
\draw[usual][postaction={decorate,decoration={markings,mark=at position 0.5 with {\node[circle,fill=red,inner sep=1.5pt] {};}}}] (1,2) [upup=1];
\draw [usual][postaction={decorate,decoration={markings,mark=at position 0.5 with {\node[circle,fill=red,inner sep=1.5pt] {};}}}] (3,0) [dnup=0];
\end{tikzpicture} 
\quad \rightarrow \quad 
\begin{tikzpicture}[anchorbase, scale=0.5]
\draw [draw=none] (1,0) rectangle (2,2);
\draw[usual][postaction={decorate,decoration={markings,mark=at position 0.5 with {\node[circle,fill=red,inner sep=1.5pt] {};}}}] (1,2) [upup=1];
\draw[usual][postaction={decorate,decoration={markings,mark=at position 0.5 with {\node[circle,fill=red,inner sep=1.5pt] {};}}}] (3,2) [upup=1];
\end{tikzpicture} 
.
\end{gather*}
As for Temperley--Lieb, we collect the resulting diagrams into what we call \emph{walk basis}. 

\begin{Theorem}\label{T:Basis}
Interpreting the blue dot as in \autoref{S:Dia} and the red dot as the projector in \autoref{S:Finite} (details in the proof) and cups and caps as usual,
the walk basis gives a basis of $\text{End}_G(V^{\otimes k})$.
\end{Theorem}

\begin{proof}
Let us first discuss the case of only blue dots. Note that the procedure describe above implies that there is an even number of blue dots, two blue dots are never on the same strand and all blue dots are reachable from the left. Thus, by \autoref{C:EqInf} and the results leading to it, we can associate blue dotted cup-caps to the endomorphism in $\text{End}_G(V^{\otimes k})$ as discussed in the preparation of \autoref{C:EqInf}. So, since the dimension of $\text{End}_G(V^{\otimes k})$ is $\frac{1}{2}\binom{2k}{k}$ by \autoref{lemma:dimension-dicinf} as is the number of closed walks of length $2k$ on the infinite type $D$ Dynkin diagram (by \autoref{L:Count} below), it remains to argue that every diagram predicted by \autoref{L:tlbasis} appears exactly once. This however follows from the Temperley--Lieb case by redirecting the walks created for it to the vertex $0'$ whenever possible.

Now we turn our attention to the red dots, where the same strategy as for blue dots applies. Precisely, the procedure ensures that there is an even number of red dots, two red dots are never on the same strand and all red dots have $n-2$ through strands in between them and the far left. For example,
\begin{gather*}
\tilde{D}_{3}=
\dynkin[
labels={0',0,1,2,3',3},
label directions={,,left,right,,},
scale=1.8,
extended
] D{ooo*o},\quad
w=(+,+,\rplus,\rminus,-,-)
\leftrightsquigarrow
\begin{tikzpicture}[anchorbase, scale=0.5]
\draw [draw=none] (1,0) rectangle (2,2);
\draw [usual] (0,0) [dnup=0];
\draw [usual][postaction={decorate,decoration={markings,mark=at position 0.5 with {\node[circle,fill=red,inner sep=1.5pt] {};}}}](1,0) [dndn=1];
\draw [usual][postaction={decorate,decoration={markings,mark=at position 0.5 with {\node[circle,fill=red,inner sep=1.5pt] {};}}}](1,2) [upup=1];
\end{tikzpicture}
.
\end{gather*}
This means we can replace every dotted cup-cap by a projector box, after appropriate Morse positioning (the relations imply that the choice of such a positioning will not affect the endomorphism we want to associate to the diagram). For example, for $n=3$,
\begin{gather*}
\begin{tikzpicture}[anchorbase, scale=0.5]
\draw [draw=none] (1,0) rectangle (2,2);
\draw [usual] (0,0) [dnup=0];
\draw [usual][postaction={decorate,decoration={markings,mark=at position 0.5 with {\node[circle,fill=red,inner sep=1.5pt] {};}}}](1,0) [dndn=1];
\draw [usual][postaction={decorate,decoration={markings,mark=at position 0.5 with {\node[circle,fill=red,inner sep=1.5pt] {};}}}](1,2) [upup=1];
\end{tikzpicture}
\leftrightsquigarrow
\begin{tikzpicture}[anchorbase,scale=1]
\draw[ipp] (0,0) rectangle (1,0.5) node[pos=0.5]{$J$};
\end{tikzpicture}
.
\end{gather*}
Using the results of \autoref{S:Finite}, this gives us an interpretation of dotted diagrams as maps in $\text{End}_G(V^{\otimes k})$.
Now, the dimension count comes out as expected by \autoref{T:CharFormula} and \autoref{L:Count}, and the rest follows as for blue dots.

\begin{Lemma}\label{L:Count}
The number of closed walks of length $2k$ on the infinite type $D$ Dynkin diagram
(starting and ending at a fixed leaf of the fork) is $\frac{1}{2}\binom{2k}{k}$.

Moreover, the number of closed walks of length $2k$ on the affine type $D$ Dynkin diagram
(starting and ending at a fixed leaf of the fork) is, for $k\ge1$,
\begin{gather*}
\frac{2^{2k+1}}{4n}\left( 1 + \sum\limits_{m=1}^{n-1} \cos^{2k}\left(\frac{m\pi}{n}\right)\right).
\end{gather*}
\end{Lemma}

\begin{proof}
This well-known count is surprisingly difficult to find in the literature, so we give a quick proof. The two parts (infinite and affine) are proven using two different (yet related) methods: via a generating function and via a spectral decomposition.

\textit{Infinite case.}
Count in pairs of steps and set $P(x)=\sum_{k\ge0}p_k x^k$, where $p_k$ is the
number of closed walks of length $2k$ based at $c$.
Let $T(x)$ be the even return series at the first tail vertex $t_1$, avoiding the
edge back to $c$. An excursion is empty or goes one step away, makes an excursion,
and returns, hence
\[
T(x)=1+xT(x)^2 \qquad\Rightarrow\qquad T(x)=\frac{1-\sqrt{1-4x}}{2x}.
\]
A primitive return from $c$ chooses one of the three branches:
either one of the two leaves (contribution $x$ each), or the tail
(contribution $xT(x)$). Thus the series of primitive returns is $x(2+T(x))$,
and concatenating such blocks gives
\[
P(x)=\frac{1}{1-x(2+T(x))}
= \frac{2}{(1+\sqrt{1-4x})\sqrt{1-4x}}
= \frac{1}{2x}\Bigl((1-4x)^{-1/2}-1\Bigr).
\]
Using $(1-4x)^{-1/2}=\sum_{m\ge0}\binom{2m}{m}x^m$, we get
\[
P(x)=\sum_{k\ge0}\frac12\binom{2k+2}{k+1}x^k
=\sum_{k\ge0}\binom{2k+1}{k+1}x^k,
\]
so $p_k=\binom{2k+1}{k+1}$ as claimed.

\textit{Affine case.}
For the affine type D Dynkin diagram, let $A$ be the adjacency matrix. Decompose
\[
e_0=\frac{1}{\sqrt2}(s_0+a_0),\qquad
s_0=\frac{e_0+e_{0'}}{\sqrt2},\quad a_0=\frac{e_0-e_{0'}}{\sqrt2}.
\]
Since $A(e_0-e_{0'})=0$, we have $Aa_0=0$, hence for $k\ge1$,
\[
p_k=(A^{2k})_{00}=\langle e_0,A^{2k}e_0\rangle=\frac12\langle s_0,A^{2k}s_0\rangle.
\]

On the symmetric subspace spanned by
$s_0,e_1,\dots,e_{n-1},s_n$ with $s_n=(e_n+e_{n'})/\sqrt2$,
the operator $A$ is represented by the $(n{+}1)\times(n{+}1)$ matrix $B$ of the
weighted walk
\[
s_0 \xleftrightarrow{\ \sqrt2\ } 1-2-\cdots-(n-1)\xleftrightarrow{\ \sqrt2\ } s_n.
\]
For $\theta=j\pi/n$ ($0\le j\le n$) define $v^{(j)}\in\mathbb{R}^{n+1}$ by
\[
v^{(j)}_0=1,\qquad v^{(j)}_m=\sqrt2\cos(m\theta)\ (1\le m\le n-1),\qquad
v^{(j)}_n=\cos(n\theta)=(-1)^j.
\]
Using $\cos((m-1)\theta)+\cos((m+1)\theta)=2\cos\theta\cos(m\theta)$ and
$1+\cos(2\theta)=2\cos^2\theta$, one checks that
\[
Bv^{(j)}=2\cos(\tfrac{j\pi}{n})\,v^{(j)}.
\]
Moreover,
\[
\|v^{(0)}\|^2=\|v^{(n)}\|^2=2n,\qquad \|v^{(j)}\|^2=n\ \ (1\le j\le n-1),
\]
so the weight of the $j$th eigenvalue at $s_0$ equals
$|v^{(j)}_0|^2/\|v^{(j)}\|^2=\frac1{2n}$ for $j=0,n$ and $\frac1n$ otherwise.
Therefore
\[
\langle s_0,B^{2k}s_0\rangle
=\frac1{2n}\,2^{2k}+\frac1n\sum_{j=1}^{n-1}\bigl(2\cos(\tfrac{j\pi}{n})\bigr)^{2k}
+\frac1{2n}\,2^{2k}.
\]
Dividing by $2$ gives the claimed formula for the number of walks.
\end{proof}

The proof is complete.
\end{proof}

By bending \autoref{L:Bending}, we also obtain from \autoref{T:Basis} a basis for all hom-spaces in the dicyclic categories.

\section{Example computations in {\Magma}}

In this section we give some sample {\Magma} code for computing in the dicyclic Temperley--Lieb category.
{\Magma} is a computer algebra system designed to solve problems in algebra, number theory, geometry and combinatorics. It was developed by Cannon and their team at the University of Sydney and version 1 was first released in August 1993. We refer to the handbook for more information: \cite{BoCa-handbook}; a beginner guide from which we stole the way to display the
{\Magma} code can be found here \cite{TaTu-magma}.

\begin{Remark}\label{R:MagmaExact}
Computer algebra systems perform exact calculations; in particular, one can use {\Magma} output 
in papers or theses without losing the exactness.
\end{Remark}

\begin{Remark}
{\Magma} is a non-commercial system, but the costs (such as preparation of user documentation, the fixing of bugs, and the provision of user support)
need to be recovered. So {\Magma} is \emph{non-commercial but not free}, and the 
distribution is organized on a subscription basis. In order to get {\Magma}
on your machine, use this site: \url{http://magma.maths.usyd.edu.au/magma/ordering/}

Free, very useful, and completely sufficient 
for our calculations, is the \emph{online calculator} \url{http://magma.maths.usyd.edu.au/calc/}. All the calculations below can be copied and pasted directly into the online calculator. 
\end{Remark}

{\Magma} can compute character tables and perform character inner product calculations efficiently. For example, the code below produces the character table for the degree four dicyclic group and computes the dimensions of the endomorphism algebras $\text{End}_{G}(V^{\otimes k})$ for $1 \leq k \leq 8$:

{\color{purple}
\begin{verbatim}
>G:=SmallGroup(16,9);
>C:=CharacterTable(G);
>C;
>V:=C[6];
>for i in [1..5] do
>InnerProduct(C[6]^(2*i), C[1]);
>end for;
\end{verbatim}}

\result

\begin{verbatim}
Character Table of Group G
-------------------------------
Class |   1  2  3  4  5   6   7
Size  |   1  1  2  4  4   2   2
Order |   1  2  4  4  4   8   8
-------------------------------
p  =  2   1  1  2  2  2   3   3
-------------------------------
X.1   +   1  1  1  1  1   1   1
X.2   +   1  1  1  1 -1  -1  -1
X.3   +   1  1  1 -1  1  -1  -1
X.4   +   1  1  1 -1 -1   1   1
X.5   +   2  2 -2  0  0   0   0
X.6   -   2 -2  0  0  0  Z1 -Z1
X.7   -   2 -2  0  0  0 -Z1  Z1
1
3
10
36
136
528
2080
8256
\end{verbatim}

The printed sequence of dimensions matches that of the $n=4$ formula from \autoref{S:RepTheory}. More importantly, we can use {\Magma} to check relations between algebra generators. The code below identifies the decorated cup-cap in the endomorphism algebra (again here $n=4$) and computes the defining relations for the dicyclic Temperley--Lieb category. To speed up the calculations, we use the smallest field over which our group is semisimple.

{\color{purple}
\begin{verbatim}
>G:=SmallGroup(16,9);
>F:=GF(17);
>V:=simpleModules(G,F)[6];
//setting up the standard representation
>CUP:=Matrix(F,4,1,[0,-1,1,0]);
>CAP:=Matrix(F,1,4,[0,1,-1,0]);
>CC:=KroneckerProduct(CAP,CUP);
//defining the cup and cap morphisms as matrices
>T:=TensorPower(V,2);
>A:=EndomorphismAlgebra(T);
>P:=PrimitiveIdempotents(A);
>P;
>dCC:=-2*P[2];
//defining the decorated cup-cap appropriately
>IsZero(dCC*CUP);
>IsZero(CAP*dCC);
>IsZero(dCC^2+2*dCC);
>KroneckerProduct(ScalarMatrix(F,2,1),CAP)*KroneckerProduct(dCC,ScalarMatrix(F,2,1)) 
*KroneckerProduct(ScalarMatrix(F,2,1),CUP);
>KroneckerProduct(CAP,ScalarMatrix(F,2,1))*KroneckerProduct(ScalarMatrix(F,2,1),dCC) 
*KroneckerProduct(CUP,ScalarMatrix(F,2,1));
//confirming the relations of the dicyclic temperley-lieb category
\end{verbatim}}

The above code prints the list of primitive idempotents in $\text{End}_{G}(V^{\otimes 2})$, the first two in the list being scalar multiples of the cup-cap and decorated cup-cap (over $\mathbb{F}_{17}$ here), and the last being the second type D Jones--Wenzl projector $f_2$:

\result

\begin{verbatim}
[ 0  0  0  0]
[ 0  9  8  0]
[ 0  8  9  0]
[ 0  0  0  0],

[ 0  0  0  0]
[ 0  9  9  0]
[ 0  9  9  0]
[ 0  0  0  0],

[ 1  0  0  0]
[ 0  0  0  0]
[ 0  0  0  0]
[ 0  0  0  1]
true
true
true
[ 1  0]
[ 0  1]
[ 1  0]
[ 0  1]
\end{verbatim}

We can also verify the diagrammatic formulas for the type D Jones--Wenzl projectors. Using the same formulas for the cup, cap, and decorated cup-cap morphisms as above, we can check that the type D projectors in the dicyclic Temperley--Lieb category correspond to the correct idempotents in $\text{End}_G(V^{\otimes k})$. The code below checks the case of $f_3$:

{\color{purple}
\begin{verbatim}
>e1:=KroneckerProduct(CC,ScalarMatrix(F,2,1));
>e2:=KroneckerProduct(ScalarMatrix(F,2,1),CC);
//defining the e1,e2 generators in the 3 strand algebra
>d1:=ScalarMatrix(F,8,1);
>d2:=e2;
>d3:=(3/4)*(e1+KroneckerProduct(dCC,ScalarMatrix(F,2,1)));
>d4:=(1/2)*(e2*e1+e2*KroneckerProduct(dCC,ScalarMatrix(F,2,1))+ 
e1*e2+KroneckerProduct(dCC,ScalarMatrix(F,2,1))*e2);
>d5:=(1/4)*(e1*e2*KroneckerProduct(dCC,ScalarMatrix(F,2,1))+
KroneckerProduct(dCC,ScalarMatrix(F,2,1))*e2*e1);
>b3:=d1+d2+d3+d4+d5;
//the identity, e1, and other three terms in the expansion of b3
>b3;
\end{verbatim}}

The matrices $d1,d2,\dots, d_5$ represent the endomorphisms corresponding to each term in the diagrammatic expansion of $f_3$. The matrix printed by {\Magma} is the projection onto the submodule generated by $v_1 \otimes v_1 \otimes v_1$ and $v_2 \otimes v_2 \otimes v_2$ as we hoped:

\result

\begin{verbatim}
[ 1  0  0  0  0  0  0  0]
[ 0  0  0  0  0  0  0  0]
[ 0  0  0  0  0  0  0  0]
[ 0  0  0  0  0  0  0  0]
[ 0  0  0  0  0  0  0  0]
[ 0  0  0  0  0  0  0  0]
[ 0  0  0  0  0  0  0  0]
[ 0  0  0  0  0  0  0  1]
\end{verbatim}

As discussed previously, this submodule is the unique two dimensional simple with multiplicity one as a summand (as long as $n>3$, in this case). In general, one can produce the sequence of type D Jones--Wenzl projectors $\{f_j\}$ by applying the recursion to the matrix $\text{diag}(1,0,0,1)$ appearing on the previous page, as this matrix represents $f_2$. {\Magma} will print a $2^j \times 2^j$ matrix with 1's in the top left and bottom right corner and zeros elsewhere, as expected. 
Finally, we can also use {\Magma} to investigate the relations between higher projectors. We define the generator $J$ as the primitive idempotent corresponding to the simple representation $V_{n'}$, which is a summand of the 4th tensor power in the case of $n=4$.

{\color{purple}
\begin{verbatim}
>b2:=P[3];
>b2:=KroneckerProduct(b2,ScalarMatrix(F,2,1));
>b3:=b2+b2*e2*b2;
>b3:=KroneckerProduct(b3,ScalarMatrix(F,2,1));
>b4:=b3+b3*e3*b3;
//creates the first few generating jones-wenzl projectors
>e4:=KroneckerProduct(ScalarMatrix(F,8,1),CC);
>b:=PrimitiveIdempotents(EndomorphismAlgebra(TensorPower(V,4)))[4];
//defines the generator J, corresponding to the new summand at the nth power
>be4b:=KroneckerProduct(b,ScalarMatrix(F,2,1))*e4*KroneckerProduct(b,ScalarMatrix(F,2,1));
>b:=KroneckerProduct(b,ScalarMatrix(F,2,1));
>IsZero(-2*be4b-b);
//checks relation (7)
\end{verbatim}}

The code above corroborates the relation (7) of $\mathbf{DiTL}(\delta,n)$. All relations in the degree $n$ dicyclic Temperley--Lieb category can be checked in this manner. Below, we present code which checks all of the other relations in the degree $4$ dicyclic Temperley--Lieb category. 

{\color{purple}
\begin{verbatim}
>IsZero(KroneckerProduct(b2,ScalarMatrix(F,2,1))*b3
-b3*KroneckerProduct(b2,ScalarMatrix(F,2,1)));
>IsZero(KroneckerProduct(CC,ScalarMatrix(F,4,1))*b3);
>IsZero(KroneckerProduct(b2,ScalarMatrix(F,2,1))*e3
-e3*KroneckerProduct(b2,ScalarMatrix(F,2,1)));
>IsZero(e3*b3*e3+KroneckerProduct(b2,ScalarMatrix(F,2,1))*e3);
>b4:=PrimitiveIdempotents(EndomorphismAlgebra(TensorPower(V,4)))[5];
>b4:=KroneckerProduct(b4,ScalarMatrix(F,2,1));
>IsZero(b4+2*b4*e4*b4);
//the above are relation checks for JW projectors in the usual TL category
>IsZero(KroneckerProduct(b3,ScalarMatrix(F,2,1))*b-b);
//relation (1), j=3
>IsZero(b4*b);
//relation (3), n=4;
>IsZero(e4*b*e4+(1/2)*KroneckerProduct(b3,ScalarMatrix(F,2,1))*e4);
//relation (4), n=4
>IsZero(KroneckerProduct(e3,ScalarMatrix(F,2,1))*b);
>e5:=KroneckerProduct(ScalarMatrix(F,16,1),CC);
>IsZero(e5*KroneckerProduct(b,ScalarMatrix(F,2,1))
-KroneckerProduct(b,ScalarMatrix(F,2,1))*e5);
//relations (5),(6) i=3, i=5, n=4
\end{verbatim}}

Running this code results in {\Magma} printing one true for each \code{IsZero()} command.

\end{document}